\documentclass[11pt,english]{article}
\usepackage{hyperref}
\hypersetup{
    colorlinks=true,
    linkcolor=blue,
    filecolor=magenta,      
    citecolor=blue,
    urlcolor=black
    }
\usepackage[sc]{mathpazo}
\usepackage[T1]{fontenc}
\usepackage[latin9]{inputenc}
\usepackage{booktabs}
\usepackage{geometry}
\usepackage{color}
\usepackage{array}
\usepackage{bbding}
\usepackage{multirow}
\usepackage{algorithm}
\usepackage{algpseudocode}
\usepackage[fleqn]{amsmath}
\usepackage{exscale}
\usepackage{relsize}
\usepackage{amssymb}
\usepackage{amsthm}
\usepackage{graphicx}
\usepackage[inkscapelatex=false]{svg}
\usepackage{natbib}
\usepackage{lscape}
\usepackage{comment}
\usepackage{bm}
\usepackage[title]{appendix}
\usepackage{longtable}
\usepackage{mathtools}
\usepackage{mathrsfs}
\usepackage{chngcntr}
\usepackage{authblk}
\usepackage{tikz}
\usepackage{threeparttable}
\usepackage{makecell}

\makeatletter
\allowdisplaybreaks

\@ifundefined{showcaptionsetup}{}{%
 \PassOptionsToPackage{caption=false}{subfig}}
\usepackage{subfig}
\makeatother

\usepackage{babel}
\usepackage{enumitem}

\newtheorem{definition}{Definition}

\makeatletter
\newenvironment{breakablealgorithm}
  {
   \begin{center}
     \refstepcounter{algorithm}
     \hrule height.8pt depth0pt \kern2pt
     \renewcommand{\caption}[2][\relax]{
       {\raggedright\textbf{\ALG@name~\thealgorithm} ##2\par}%
       \ifx\relax##1\relax 
         \addcontentsline{loa}{algorithm}{\protect\numberline{\thealgorithm}##2}%
       \else 
         \addcontentsline{loa}{algorithm}{\protect\numberline{\thealgorithm}##1}%
       \fi
       \kern2pt\hrule\kern2pt
     }
  }{
     \kern2pt\hrule\relax
   \end{center}
  }
\makeatother

\begin{document}\sloppy

\title{Planning of Truck Platooning for Road-Network Capacitated Vehicle Routing Problem}

\author[a,b,c]{Yilang Hao}
\author[a,b]{Zhibin Chen \footnote{Corresponding author. zc23@nyu.edu.}}
\author[c,d]{Xiaotong Sun}
\author[e]{Lu Tong}
\affil[a]{\small\emph{Shanghai Key Laboratory of Urban Design and Urban Science, NYU Shanghai, 567 West Yangsi Road, Pudong New District, Shanghai, 200126, China}}
\affil[b]{\small\emph{Shanghai Frontiers Science Center of Artificial Intelligence and Deep Learning, NYU Shanghai, 567 West Yangsi Road, Pudong New District, Shanghai, 200126, China}}
\affil[c]{\small\emph{Department of Civil and Urban Engineering, New York University, 6 MetroTech Center, Brooklyn, NY 11201, United States}}
\affil[c]{\small\emph{Intelligent Transportation Thrust, Systems Hub, The Hong Kong University of Science and Technology (Guangzhou), Guangzhou, China}}
\affil[d]{\small\emph{Department of Civil and Environmental Engineering, The Hong Kong University of Science and Technology, Clear Water Bay, Hong Kong SAR, China}}
\affil[e]{\small\emph{The School of Traffic and Transportation, Beijing Jiaotong University, No.3 Shangyuan Cun, Haidian District, Beijing, China}}
\date{}
\maketitle

\section*{Acknowledgements}
This work was partially supported by Shanghai Chenguang Program (21CGA72), the National Natural Science Foundation of China (72201073, 72001021, 72061127003), the Preparation Fund of Shanghai Key Laboratory of Urban Design and Urban Science (Grant No. 10407\_Key Lab\_Preparation Fund), Shanghai Eastern Scholar Program (QD2020057), the NYU Shanghai Doctoral Fellowships, and the NYU Shanghai Boost Fund.

\begin{abstract}
Truck platooning, a linking technology of trucks on the highway, has gained enormous attention in recent years due to its benefits in energy and operation cost savings. However, most existing studies on truck platooning limit their focus on particular scenarios that each truck can serve only one customer demand and is thus with a specified origin-destination pair, so only routing and time schedules are taken into account. Nevertheless, in real-world logistics, each truck may need to serve multiple customers located at different places, and the operator managing a fleet of trucks thus has to determine not only the routing and time schedules of each truck but also the set of customers allocated to each truck and their sequence to visit. This is well known as a capacitated vehicle routing problem with time windows (CVRPTW), and considering the application of truck platooning in such a problem entails new modeling frameworks and tailored solution algorithms. In light of this, this study makes the first attempt to optimize the truck platooning plan for a road-network CVRPTW in a way to minimize the total operation cost, including vehicles' fixed dispatch cost and energy cost, while fulfilling all delivery demands within their time window constraints. Specifically, the operation plan will dictate the number of trucks to be dispatched, the set of customers, and the routing and time schedules for each truck. In addition, the modeling framework is constructed based on a road network instead of a traditional customer node graph to better resemble and facilitate the platooning operation. A 3-stage algorithm embedded with a "route-then-schedule" scheme, dynamic programming, and modified insertion heuristic, is developed to solve the proposed model in a timely manner. Numerical experiments are conducted to validate the proposed modeling framework, demonstrate the performance of the proposed solution algorithm, and quantify the benefit brought by the truck platooning technology.

\hfill\break%
\noindent\textit{Keywords}: Truck platooning; Vehicle routing problem; Route then schedule; Dynamic programming; Modified insertion heuristic

\end{abstract}

\section{Introduction}
In the era of advanced logistics and delivery systems, truck platooning, a concept that leverages advancements in vehicle automation and wireless communication, has emerged as a solution to revolutionize cargo delivery operations. 
Specifically, integrated with modern logistics, truck platooning is expected to yield significant benefits, including enhanced capacity utilization by closer headways \citep{chen2016, chen2017}, reduced operational costs \citep{metropolitan, talogistics}, relieved drivers' driving fatigue \citep{larsen2019hub,sun2021auction, sun2021decentralized, Janssen2015TruckPD}, and a more stable and effective supply chain \citep{supplychain}, thereby bringing substantial socioeconomic welfare and expenditure savings in logistics services. 
Specifically, the reduction in total logistics is reported to be 1\% to 2\% \citep{platoonimpact}, which is remarkable considering the huge volume of logistics. 
In addition, with the platooning function, the energy saving can be as high as more than 8\% \citep{platoonsaving, sivan2022, peloton}, though it may vary from type of vehicle, load, and platoon size.
Given the great promises held by truck platooning, companies such as Daimler AG, Volvo, DAF, TuSimple, Honda, and Peleton have already been aggressively testing this technology. For instance, TuSimple and Waymoe have, respectively, demonstrated their driverless platooning products, called TuSimple Connect and Waymo Via \citep{tusimple,waymo}. PATH, a pioneer in truck platooning and collaborated with Volve Group, presented the practical operation of truck platooning at the Port of Los Angeles, and on interstate highway I-66 in the northern Virginia suburbs of Washington DC in 2017 \citep{path}. European Automobile Manufacturer's Association also provided a blueprint describing how to achieve multi-brand platooning before 2025. 

In addition to the industry, many scholars have been studying the truck platooning operation, and we will first review them to identify the research gaps that motivate this study. Given that forming a platoon entails trucks arriving at the same location at the same time, both routing and time scheduling plans matter. Regarding the routing plan, \cite{LARSSON2015258} formulated a graph routing problem based on truck platooning but ignored the constraints of delivery deadlines to reduce the problem's complexity. Two constructive heuristics and a local search algorithm are proposed to optimize each truck's route for forming platoons in a way to minimize the system cost. Based on a similar logic, \cite{Larson2013} further optimized the velocity of vehicles to maximize fuel savings by platooning, and derived a hub-based heuristic to simplify the solving process. Among the studies concerning the time schedule, \cite{boysen201826} optimized trucks' departure times to form platoons based on a simple scenario in that all trucks adopt an identical path. Therefore, the optimal solution to minimize operation cost can be found in a timely manner, and they further explored that tight time windows and limited platoon size may weaken the platooning benefit. \cite{Hoef2018} studied a similar problem but allowed speed controls for the rear vehicle, and they proposed a local improvement heuristic to tackle large-scale instances. As an extension to \cite{Hoef2018}, \cite{hoef20174228} constructed a stochastic dynamic programming model to quantify the influence of time uncertainty. \cite{zhang20171} is another study that considers platooning with time uncertainties, and they discovered that platooning is only preferable when the travel time interval of two vehicles is within a particular threshold and that the uncertainty of travel time further reduces this threshold. 

In order to facilitate the truck platooning in a network, nevertheless, it is necessary to jointly consider the routing and time scheduling plans of each truck. In this regard, \cite{Larson2016CoordinatedPR} constructed a mixed integer program to optimize the routing and time scheduling plans of trucks, and applied a model decomposition approach to seize the exact solution. In addition, a preprocessing process was developed to eliminate unnecessary paths and thus reduce the problem size. \cite{LUO2018213} extended the above study by considering the adjustment of trucks' speeds. \cite{mandatorybreak} considered mandatory breaks of drivers during travel of a platooning problem to capture the real-world scenario better, and a partial mixed integer program approach integrated with the iterated neighborhood search method was developed to solve the problem efficiently. \cite{Luo_2022} innovatively invented a "route-then-schedule" iterative process to sequentially solve the routing problem and time scheduling problem separately, thus reducing the problem complexity during each solving process. They also provided analytical proof that such an iterative process will converge to a suboptimal solution with their designed feedback link cost modification mechanism. Recently, \cite{zhao2023improved} resorted to the above "route-then-schedule" solution framework to solve a truck routing problem considering that leading trucks also have fuel savings and trucks are not allowed to wait during the trip. Instead of setting time consistency constraints to record schedules of trucks, \cite{ABDOLMALEKI202191} set up a time-expanded network to schedule the travel itineraries of trucks to form platoons for saving energy, and the problem was formulated as a multicommodity flow problem. Approximation algorithms together with dynamic programming were then proposed to seize a decent solution.

For the aforementioned studies on truck platooning, each truck is assumed to serve only one predetermined customer demand, and thus it has only one origin and destination which are both predetermined as well. Nevertheless, in the real-world logistics or delivery service, each truck is capable of serving multiple customers along its itinerary. In addition, the set of customers to be served by each truck and their sequence to receive service are yet to be decided. In this study, we make the first attempt to optimize the truck platooning plan for the real-world delivery problem. Particularly, without loss of generality, we assume that there are given customers located at different places on the network. Each customer has a specific demand and a required time window to be served, and each truck has a maximum capacity of goods to carry. Acting as a central coordinator or a carrier company, our goal is to optimize the truck platooning plan in a road network to serve all the customer demands in a way that minimizes the total cost, including the dispatch cost and energy cost. Specifically, the plan will dictate not only the routing and time schedules of each truck but also the set of customers allocated to each truck and their sequence to serve. Unlike traditional vehicle routing problems (VRPs), the problem considered in this study is modeled on the road network instead of the customer node graph so as to delineate the platoon formation and disassembling. In addition, as the weight of a truck will change during its trip and affect the energy saving brought by the truck platooning, its impact on the energy consumption rate is considered explicitly. Therefore, this study essentially focuses on the planning of truck platooning for a road-network capacitated vehicle routing problem with time windows (RNCVRPTW). 

\begin{figure}[htbp]
	\centering
    \includegraphics[scale=0.7]{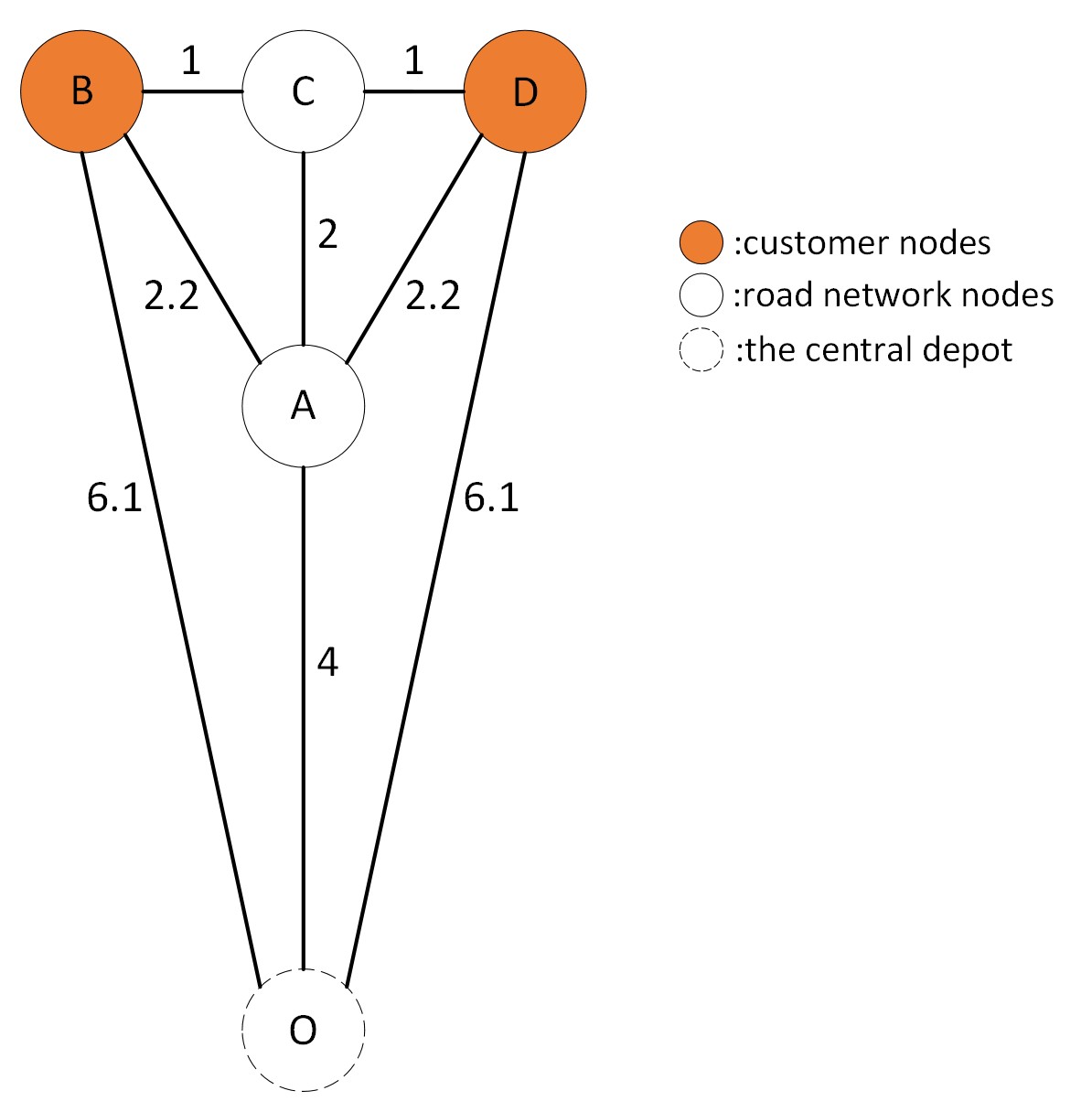}
	\caption{Toy network}
    \label{fig:toy}
\end{figure}

\begin{figure}[htbp]
	\centering
	\subfloat[Without platooning]{\label{fig:illusb}\includegraphics[scale=0.7]{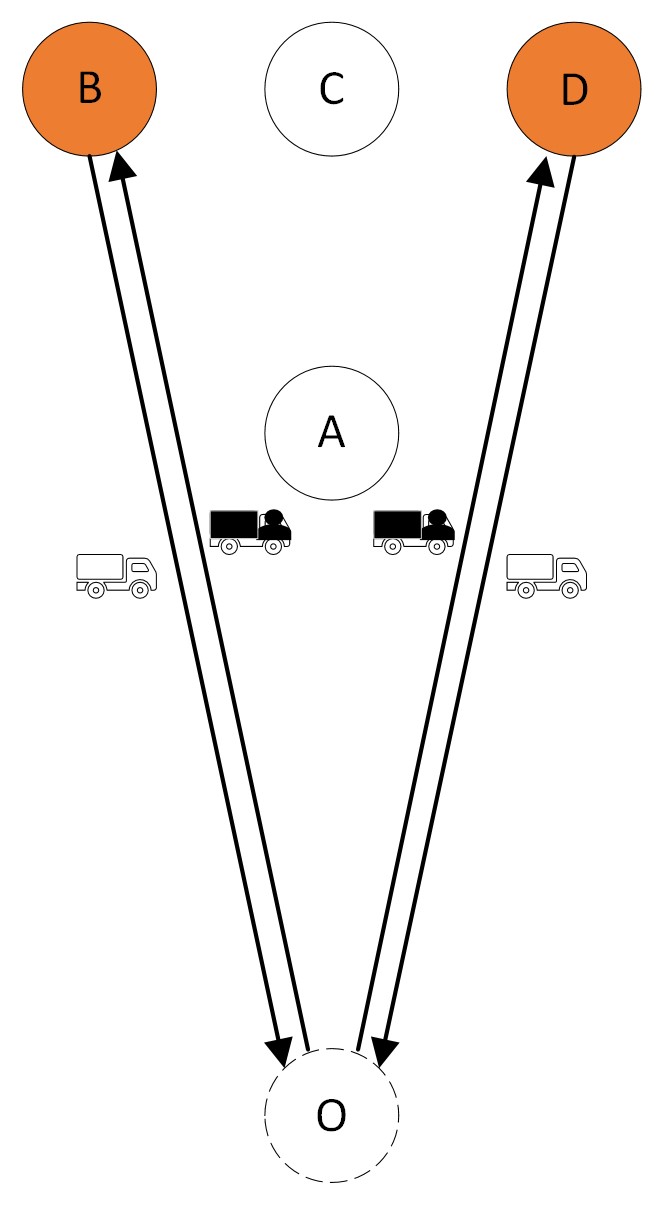}}\quad
    \subfloat[With platooning]{\label{fig:illusa}\includegraphics[scale=0.7]{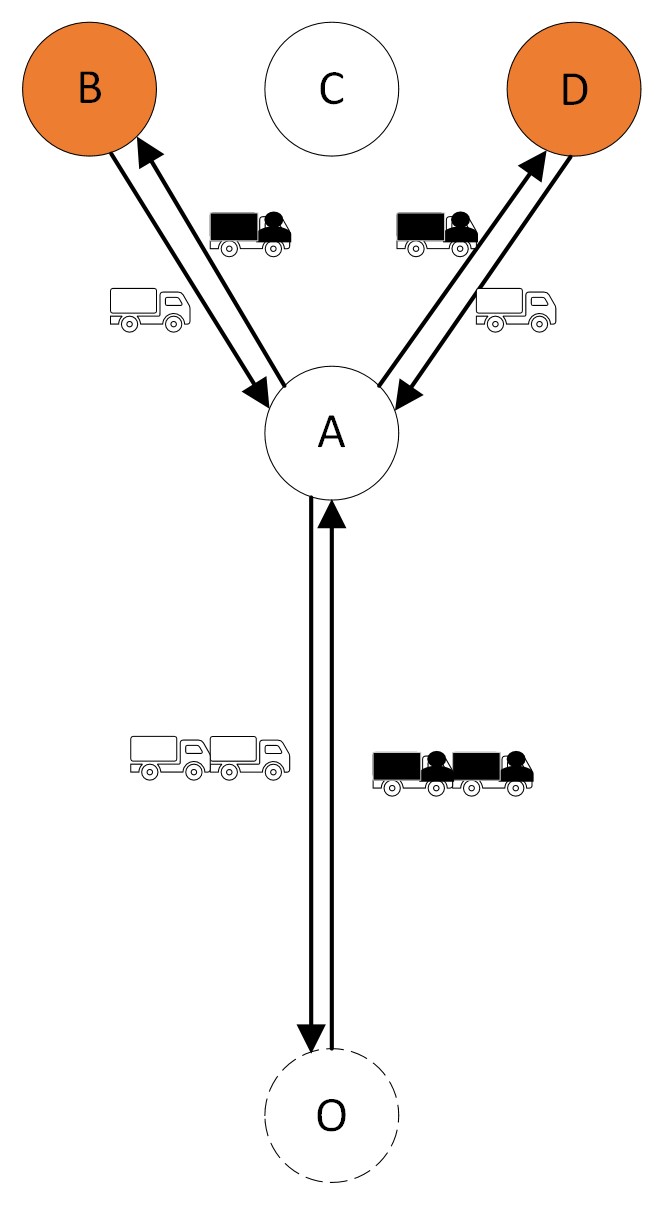}}\\
	\caption{Optimal trajectories with or without platooning}
    \label{fig::illus}
\end{figure}

Below a toy example is presented to demonstrate the impact of truck platooning on the delivery problem. The network, shown in Figure \ref{fig:toy}, is composed of 5 vertices, with 1 depot at node $O$, 2 colored customer nodes which are $B$ and $D$, and 2 intermediate nodes which are $A$ and $C$. The link travel times are provided along their corresponding links, and for simplicity, we assume that the fuel consumption rate with respect to time is 1, and the demand in each customer node is exactly equal to the truck capacity. As a result, two trucks are needed to fulfill the demands. In this example, we omit time windows for simplicity, and the energy-saving ratio for the following trucks in a platoon is chosen as 10\%, as \cite{peloton} and \cite{bishop2020} reported. We also assume that a full-load truck consumes 20\% more fuel than an empty truck, as indicated by \cite{weighteffect}, \cite{weighteconomy}, and \cite{quantitativeeffect}. Without loss of generality, we set the fuel consumption rate as 1 per unit travel time for an empty truck, so the one for a full-load truck is 1.2.  Given the above setting, the optimal routes of dispatched trucks with and without the support of platooning technologies are visually presented in Figure \ref{fig::illus}. As shown in Figure \ref{fig:illusb}, under the no-platoon case, two trucks are dispatched to serve two customers, respectively, yielding a total energy cost of 
$$(6.1\times 1.2 + 6.1\times 1)\times 2 = 26.84$$
where $6.1\times 1.2$ represents the fuel consumption on $(O,B)$ or $(O,D)$ when the truck is fully-loaded, and $6.1\times 1$ represents the one on $(B,O)$ or $(D,O)$ when the truck is empty. 

In Figure \ref{fig:illusa}, as the platooning feature is enabled, two trucks will platoon together on the common roads $(O,A)$ and $(A,O)$ to save energy. Accordingly, the total fuel cost considering platooning is:
$$4\times (1+1-10\%)\times (1+1.2) + 2.2\times (1+1.2)\times 2 = 26.4$$
where $(1+1-10\%)$ appeared in the first term takes the platoon savings on link $(O,A)$ and $(A,O)$ into consideration. Therefore, the platooning induces a $\dfrac{26.84-26.4}{26.84} = 1.64\%$ in the energy cost. Such a simple instance illustrates that, if platooning is allowed in the truck delivery operation, the itineraries of trucks may be altered to induce additional cost savings. In addition, as highlighted by \cite{platoonseq}, the impulse of getting platoon savings can alter not only trucks' routes and schedules but also the set of customers to be served by each truck and their sequence to receive services. 

In this study, a mixed integer programming (MIP) framework is first developed to model the operation of truck platooning for RNCVRPTW. Given that the problem is NP-hard and road-network-based \citep{YAO202121,rvrp}, a 3-stage algorithm integrating the "route-then-schedule" scheme \citep{Luo_2022}, dynamic programming, and modified insertion heuristic, is developed to solve the proposed model in a timely manner with high accuracy. The extant body of literature reveals that certain local container drayage problems bear a resemblance to the issues addressed in our research. Specifically, these drayage problems delved into the delivery and relocation of containers, optimizing the routing of trucks to meet all demands while minimizing costs. A subset of these studies, as referenced in \cite{doubletrailer} and \cite{multitrailer, multitrip}, allowed for a single truck to transport two or more containers during its journey. Their problem contexts were analogous to ours in serving multiple customer demands during one trip as each container can be conceptualized as a distinct demand. However, they did not consider time-window constraints and assumed the same demand sizes for different customers. In addition, these studies related to multi-trailer drayage problems did not take the platooning feature into consideration, marking a big difference in the operation planning. There also exist several other drayage studies that discuss the operation planning under the truck platooning context, such as \cite{drayage2020}, \cite{localcontainer2021}, and \cite{drayagexie}. 
Nevertheless, these studies assumed that one truck can carry only one container or customer demand, so the operation planning of container deliveries was not fully investigated. In a nutshell, though the above container-drayage studies share some similarities with ours, their methodologies cannot be applied to our problem which optimizes not only the itineraries of trucks but also the allocations of customers to each truck and their sequence to serve.

Our contribution is three-fold: 1) This is the first study to optimize the planning of truck platooning for a road-network capacitated vehicle routing problem. It is a complicated problem because we not only take into account the road network and the potential for platooning but also the fact that there are no predetermined paths for the trucks and predetermined assignments of which customers each truck will serve. 2) We not only consider the vehicle dispatch cost and fuel consumption on the road to comprehensively depict the operational costs of a carrier company but also account for the impact of truck weight. The impact extends beyond just an increase in the rate of fuel consumption; it's also crucial to emphasize that, as energy savings are calculated on a percentage basis, the magnitude of platoon savings becomes significantly greater with an increase in weight. This provides a more accurate representation of fuel costs in the context of transportation problems. 3) We propose a novel iterative solution algorithm that combines the "route-then-schedule" scheme, dynamic programming, and modified insertion heuristic, in a way that practical large-scale problems can be solved efficiently and accurately. 

The remainder of this paper is organized as follows. The problem is mathematically formulated in Section \ref{sec:formulation}. A novel iterative solution algorithm is developed in Section \ref{sec:alg}. In Section \ref{sec:exp}, numerical experiments are conducted to demonstrate the performance of our proposed modeling framework and solution algorithm, as well as the influence of several key parameters on the system cost. Lastly, Section \ref{sec:con} concludes the paper.

\section{Problem Formulation} \label{sec:formulation}
\subsection{Problem setting}\label{sec:descrip}
In this study, delivery demands are located at a subset of nodes in the road network, and each demand is required to be delivered during the required time window, composed of the earliest arrival time and the latest departure time. In addition, each demand may have a different demand size, but the size is assumed to not exceed the capacity of each truck, and each demand cannot be served by multiple trucks to reduce the model complexity. 

We denote the set of all nodes as $N$, and it is composed of two parts, the depot $D$, and the set of all road-network nodes $R$. In this paper, only one depot is considered and we fix its index as $0$. Each truck dispatched by the carrier company must depart from the depot and return to it after serving all the demands, and each truck can only be dispatched once. The latter assumption can be relaxed by treating a truck with multiple dispatched orders as multiple trucks, and doing so imposes no changes on the modeling framework proposed in this paper. For the network nodes, they can be further divided into two categories, the customer nodes (defined by $C$), and the other intermediate nodes that can be served for platoon formation and disintegration. For each customer node, $\forall i\in C$, its demand information includes the size of demand, $q_i$, and the time window, indicating the earliest arrival time, $T^{EA}_i$, and latest departure time, $T^{LD}_i$. The set of all arcs is denoted as $A$, so our road network graph can be represented as $G=(N,A)$. Without loss of generality, we assume that all trucks travel at the same speed, and we define a set $K$ for all trucks, whose size can be chosen as the same size as $C$ to ensure the problem's feasibility. 

As mentioned before, to precisely capture the fuel consumption of trucks during their trips, we have to take the influence of the weight of trucks into consideration. \cite{linearweight} provided a substantial indication that the increase in fuel consumption rate is linearly correlated with the load of trucks. Let $\gamma$ denote the empty static weight of each truck; $\alpha$ define the corresponding fuel consumption per unit weight of an empty truck; and $\eta$ be the increased fuel consumption due to one unit of additional load filled in. According to reports and studies on truck fuel consumptions \citep{weighteffect, weighteconomy, quantitativeeffect}, the fuel consumption rate at a full-load state, which equals the truck capacity $Q$, can be 15\%-25\% larger than the basic rate. Due to the aerodynamic benefit, the following trucks in a platoon can save energy, and we define the cost-saving rate by $\beta$, which may vary from 0.05 to 0.15 in previous studies \citep{Larson2016CoordinatedPR}. There are restrictions on the platoon size \citep{laborcost} in real-world operations due to many factors such as safety concerns and the capability of communication technologies, and we denote the maximum platoon size by $L$.

\subsection{Mathematical formulation of the RCVRPTW-TP model}\label{sec:form}

Before presenting the overall model, the notations, including parameters and decision variables, frequently used in the mathematical model are provided in Table \ref{list_of_notations}.

\begin{table}[htbp]
\caption{List of Notations}
\label{list_of_notations}
\centering
\begin{tabular}{cl}
\toprule
Notation                               & \multicolumn{1}{c}{Definition}\\ 
\midrule
\multicolumn{2}{c}{\bf Sets}\\ \hline

$C$ & Set of customer nodes, $C = \left\{1,2,3,\dots, c\right\}$\\
$D$ & Depot, $D = \left\{0\right\}$\\
$R$ & Set of road-network nodes, $C \subset R$\\    
$N$ & Set of all nodes, $N = R\cup D$\\
$A$ & Set of all arcs\\
$G$ & The entire graph of the road network, $G = (N,A)$\\
$K$ & Set of trucks\\
$\mathcal{R}_k$ & The set of edges along the route adopted by truck $k\in K$ \\
$\mathcal{K}_{i,j}$ & Set of trucks taking the edge $(i,j)\in \cup_{k\in K} \mathcal{R}_k$ \\
\hline
\multicolumn{2}{c}{\bf Parameters}\\ \hline
$L$ & Maximum platoon size\\
$t_{i,j}$ & Time to travel arc $(i,j)\in A$\\
$T^{EA}_i$ & Earliest arrival time for customer node $i,i\in C$\\
$T^{LD}_i$ & Latest departure time for customer node $i,i\in C$\\
$q_i$ & Demand at customer node $i,i\in C$, and $q_i = 0$ for $i \notin C$\\
$Q$ & Truck capacity\\
$\gamma$ & Static weight of a truck\\
$\alpha$ & Base fuel consumption rate for a unit of travel time (for an empty truck)\\
$\eta$ & Marginal fuel consumption coefficient for a unit of additional weight (unit cargo weight)\\
$\beta$ & The platoon saving factor\\
$c_1$ & Vehicle fixed dispatch cost per truck\\
$c_2$ & Weight coefficient to balance two kinds of costs\\
\hline
\multicolumn{2}{c}{\bf Variables}\\ \hline
$x_{i, j, k}$ & Binary variables to indicate whether truck $k$ travels on edge $(i,j)$\\
$g_{i,k}$ & Binary variables to indicate if customer node $i$ is served by truck $k$.\\$l_{i, j, k}$ & Binary variable: 1 if truck $k$ is the leading vehicle of a platoon on edge $(i,j)$; and 0, \\
& otherwise. \\
$f_{i, j, k_1, k_2}$ & Binary variable: 1 if truck $k_1$ is the leading vehicle and $k_2$ is in the same platoon on \\
& edge $(i,j)$; and 0, otherwise.\\
$y_{i, k}$ & Volume of truck $k$ at vertex $i$\\$s_{i, k}$ & The time for truck $k$ arrives at node $i$\\
$w_{i, k}$ & The time for truck $k$ waits at node $i$\\
\bottomrule
\end{tabular}
\end{table}

\newpage
\subsubsection{Objective function}

\begin{flalign}
\min z = &c_1 \sum_{i\in D} \sum_{j\in N,(i,j)\in A} \sum_{k\in K}x_{i, j, k} +c_2 \sum_{i\in N}\sum_{j\in N,(i,j)\in A} \dfrac{\alpha}{\gamma} t_{i, j}\bigg[\sum_{k\in K}l_{i, j, k}\cdot (\eta y_{i, k} + \gamma) &\notag\\
&+(1 - \beta)(\sum_{k_1\in K}\sum_{k_2\in K:k_2 \neq k_1} f_{i, j, k_1, k_2}\cdot (\eta y_{i, k_2} + \gamma))\bigg]   &\tag{1}\label{obj}\end{flalign}
\noindent where $x_{i,j,k}$ is a binary variable to specify whether a truck $k$ will traverse edge $(i,j)$; $l_{i,j,k}$ tells if truck $k$ is a leading vehicle of a platoon on an edge $(i,j)$; $y_{i,k}$ is the volume carried by truck $k$ when it arrives at node $i$; and $f_{i,j,k_1,k_2}$ represents whether truck $k_2$ is following a leading truck $k_1$ in the same platoon on edge $(i,j)$. Obviously, any truck can only act as either a leading truck or a following truck at any time, and there is only one leading truck for any platoon. 

The above objective function minimizes the total operation cost of the carrier company, including vehicles' dispatch cost and energy costs. Specifically, the first term is the vehicles' dispatch cost, in which $\sum_{i\in D} \sum_{j\in N} \sum_{k\in K}x_{i, j, k}$ indicates the number of trucks that depart from the depot, and $c_1$ is the unit truck dispatch cost. The second term is the energy consumption cost, and it is proportional to the on-road travel time $t_{i,j}$ for any specific edge $(i,j)$. 
Specifically, for a leading truck $k$ in a platoon that travels on edge $(i,j)$, its energy consumption is calculated as $\dfrac{\alpha}{\gamma} \sum_{k\in K}l_{i, j, k}\cdot (\eta y_{i, k} + \gamma)$. On the other hand, for a following truck $k_2$ in a platoon, it can enjoy the platooning benefit, so a platoon saving factor, $1-\beta$, will be applied to obtain its actual fuel consumption, which is $\dfrac{\alpha}{\gamma}(1 - \beta) t_{i,j}(\sum_{k_1\in K}\sum_{k_2\in K:k_2 \neq k_1} f_{i, j, k_1, k_2}\cdot (\eta y_{i, k_2} + \gamma))$. To be mentioned, the factor $c_2$ is used to balance the dispatch cost and energy consumption costs in objective value. 

\subsubsection{Flow conservation constraints}
\begin{flalign*}
&\sum_{i\in N, (i,j)\in A} x_{i, j, k}=\sum_{h\in N, (j,h)\in A} x_{j, h, k} \leq 1,
&\forall j\in N, k\in K, & \tag{2}\label{2}\\
&\sum_{k\in K} g_{i, k} = 1
&\forall i\in C, &\tag{3}\label{3}\\
&\sum_{i \in N, (i,j)\in A} x_{i, j, k} \geq  g_{j, k}
&\forall j\in C,k\in K, &\tag{4}\label{4}
\end{flalign*}
where $g_{i,k}$ indicates if customer node $i$ is served by truck $k$.

Constraints (\ref{2})-(\ref{4}) construct trucks' routes and ensure that all customers are served. Constraints (\ref{2}) guarantee that the inflow of trucks should be equal to the outflow of trucks at any node, and each truck can only be dispatched at most once. Constraints (\ref{3}) and (\ref{4}) ensure that each customer node is served by exactly one truck that visits the node.

\subsubsection{Volume constraints}
\begin{flalign*}
&y_{0, k} = \sum_{i\in C} g_{i, k} \cdot q_i
&\forall k\in K, &\tag{5}\label{5}\\
&y_{0, k} \leq Q
&\forall k\in K, &\tag{6}\label{6}\\
&y_{i, k} - q_j \cdot g_{j, k} + M(1-x_{i,j,k}) \geq y_{j, k}
&\forall k\in K,  (i, j)\in A, j\in C, &\tag{7}\label{7}
\end{flalign*}
where $Q$ is the maximum capacity, and $q_j$ is the customer demand at node $j$. $M$ is a sufficiently large number. 

Constraints (\ref{5})-(\ref{7}) honor the volume consistency of each truck along its itinerary. Specifically, constraints (\ref{5}) record the cargo of each truck that needs to carry before it departs from the depot, which is fixed at node $0$ in our model. Constraints (\ref{6}) demonstrate that the initial load of each truck should never exceed its capacity. In addition, we only consider deliveries but not pickups in our problem, so the load at all following nodes will be not greater than the initial load, thereby not exceeding the capacity. Constraints (\ref{7}) specify the volume changes of each truck along its route. When it finishes delivery at a customer node, the volume of the truck will drop by the amount which is equal to the customer demand $q_j$; otherwise, the volume will remain the same.

\subsubsection{Platoon constraints}
\begin{flalign*}
&x_{i,j,k} = l_{i,j,k} + \sum_{k_1\in K,k_1\neq k}f_{i,j,k_1,k} 
&\forall (i,j)\in A, k\in K, &\tag{8}\label{8}\\
&\sum_{k_2\in  K:k_2 \neq k_1} f_{i, j, k_1, k_2} \leq (L - 1) \times l_{i, j, k_1}
&\forall (i, j)\in A, k_1\in K, &\tag{9}\label{9}\\
&-M(1 - f_{i, j, k_1, k_2}) \leq s_{i, k_1} + w_{i, k_1} - s_{i, k_2} - w_{i, k_2}
&\forall (i,j)\in A, k_1,k_2\in K, k_1\neq k_2, &\tag{10}\label{10}\\
&s_{i, k_1} + w_{i, k_1} - s_{i, k_2} - w_{i, k_2} \leq M(1 - f_{i, j, k_1, k_2})
&\forall (i,j)\in A,  k_1\neq k_2\in K, &\tag{11}\label{11}\\
&s_{i, k} + t_{i, j} + w_{i, k} - M(1 - x_{i, j, k}) \leq s_{j, k}
&\forall (i, j)\in A, j\notin D, k\in K, &\tag{12}\label{12}
\end{flalign*}
where $L$ is the maximum platoon size; $s_{i,k}$ denotes the time that truck $k$ arrives at node $i$; and $w_{i,k}$ records the time that truck $k$ waits at node $i$. Therefore, $s_{i,k}+w_{i,k}$ can be used to indicate the time that truck $k$ departs from node $i$.

Constraints (\ref{8})-(\ref{12}) honor the platooning consistency among trucks on the road. Constraints (\ref{8}) illustrate that each truck can only act as a following vehicle or a leading vehicle at the same time. Constraints (\ref{9}) specify that the length of any platoon should never exceed the size limit. Constraints (\ref{10}) and (\ref{11}) ensure that for two trucks in the same platoon, they should enter the edge $(i,j)$ at the same time. Constraints (\ref{12}) act as time consistency constraints for each truck during its trip.

\subsubsection{Time window constraints}
\begin{flalign*}
&T^{EA}_i \cdot g_{i, k} \leq s_{i, k}
&\forall i\in C, k\in K, &\tag{13}\label{13}\\
&s_{i, k} + w_{i, k} \leq T^{LD}_i \cdot g_{i, k}
&\forall i\in C, k\in K, &\tag{14}\label{14}
\end{flalign*}

Constraints (\ref{13}) and (\ref{14}) guarantee that each customer node should be served during its allowed time window, i.e., between its earliest arrival time and latest departure time.

\subsubsection{Integrity constraints}
\begin{flalign*}
&y_{i, k}\geq 0
&\forall i\in N, k\in K, &\tag{15}\label{15}\\
&x_{i,j,k}\in \left\{0,1\right\}
&\forall (i,j)\in A,k\in K, &\tag{16}\label{16}\\
&l_{i, j,k}\in \left\{0,1\right\}
&\forall (i,j)\in A,k\in K, &\tag{17}\label{17}\\
&f_{i,j,k_1,k_2}\in \left\{0,1\right\}
&\forall (i,j)\in A,k_1,k_2\in K,k_1\ne k_2 &\tag{18}\label{18}\\
&s_{i, k}\geq 0
&\forall i\in N, k\in K, &\tag{19}\label{19}\\
&w_{i, k}\geq 0
&\forall i\in N, k\in K, &\tag{20}\label{20}\\
&g_{i,k}\in \left\{0,1\right\}
&\forall i\in C,k\in K, &\tag{21}\label{21}\\
\end{flalign*}

Constraints (\ref{15}) - (\ref{21}) specify the ranges of all decision variables.

\subsubsection{Overall model}
Therefore, the road-network capacitated vehicle routing problem considering truck platooning (RCVRPTW-TP) is summarized as follows:
\begin{flalign}
&\textbf{Minimize}& &\notag \\
&Obj = c_1 \sum_{i\in D} \sum_{j\in N} \sum_{k\in K}x_{i, j, k} +c_2 \sum_{i\in N}\sum_{j\in N,j\neq i} \dfrac{\alpha}{\gamma} t_{i, j}\bigg[\sum_{k\in K}l_{i, j, k}\cdot (\eta y_{i, k} + \gamma) &\notag \\
&+(1 - \beta)(\sum_{k_1\in K}\sum_{k_2\in K:k_2 \neq k_1} f_{i, j, k_1, k_2}\cdot (\eta y_{i, k_2} + \gamma))\bigg]  &\notag\\   
&\mbox{s.t.}\quad (\ref{2})-(\ref{21})&\notag
\end{flalign}

However, the mathematical program above is a mixed integer nonlinear program because $l_{i, j, k}\cdot (\eta y_{i, k}+ \gamma)$ and $\sum_{k_1\in K}\sum_{k_2\in K:k_2 \neq k_1} f_{i, j, k_1, k_2}\cdot (\eta y_{i, k_2} + \gamma)$ in our objective function are nonlinear. To linearize them, we introduce variables $v^{(l)}_{i,j,k}$ and $v^{(f)}_{i,j,k}$ and the following constraints (\ref{24})-(\ref{29}). Specifically, constraints (\ref{22}) and (\ref{23}) define the range of variables $v^{(l)}_{i,j,k}$ and $v^{(f)}_{i,j,k}$. Constraints (\ref{24})-(\ref{26}) together ensure that $v^{(l)}_{i,j,k} = l_{i, j, k}\cdot (\eta y_{i, k}+ \gamma)$ under all circumstances. Similarly, constraints (\ref{27})-(\ref{29}) ensure that $v^{(f)}_{i,j,k} = \sum_{k_1\in K}\sum_{k_2\in K:k_2 \neq k_1} f_{i, j, k_1, k_2}\cdot (\eta y_{i, k_2} + \gamma)$ for all cases. 

\begin{flalign*}
&v^{(l)}_{i,j,k} \geq 0
&\forall (i,j)\in A, k\in K, &\tag{22}\label{22}\\
&v^{(f)}_{i,j,k} \geq 0
&\forall (i,j)\in A, k\in K, &\tag{23}\label{23}\\
&v^{(l)}_{i,j,k} \leq (\eta Q+\gamma)l_{i,j,k}
&\forall (i,j)\in A,k\in K, &\tag{24}\label{24}\\
&v^{(l)}_{i,j,k} \leq \eta y_{i,k} + \gamma 
&\forall (i,j)\in A,k\in K, &\tag{25}\label{25}\\
&v^{(l)}_{i,j,k} \geq \eta y_{i,k} + \gamma - (\eta Q+\gamma)(1-l_{i,j,k})
&\forall (i,j)\in A,k\in K, &\tag{26}\label{26}\\
&v^{(f)}_{i,j,k} \leq (\eta Q+\gamma)\sum_{k_1:k_1\neq k}f_{i,j,k_1,k}
&\forall (i,j)\in A,k\in K, &\tag{27}\label{27}\\
&v^{(f)}_{i,j,k} \leq \eta y_{i,k} + \gamma 
&\forall (i,j)\in A,k\in K, &\tag{28}\label{28}\\
&v^{(f)}_{i,j,k} \geq \eta y_{i,k} + \gamma - (\eta Q+\gamma)(1-\sum_{k_1:k_1\neq k}f_{i,j,k_1,k})
&\forall (i,j)\in A,k\in K, &\tag{29} \label{29}\\
\end{flalign*}

Accordingly, the RCVRPTW-TP model can be reformulated into a mixed-integer linear program (MILP):

\noindent \textbf{RCVRPTW-TP-R}
\begin{flalign}
&\textbf{Minimize} \notag \\
&Obj = c_1 \sum_{i\in D} \sum_{j\in N} \sum_{k\in K}x_{i, j, k} +c_2 \sum_{i\in N}\sum_{j\in N,j\neq i} \dfrac{\alpha}{\gamma} t_{i, j}[\sum_{k\in K}v^{(l)}_{i, j, k} +(1 - \beta)\sum_{k\in K} v^{(f)}_{i, j, k}]   \notag\\   
&\mbox{s.t.}\quad (\ref{2})-(\ref{29})&\notag
\end{flalign}

\section{Solution Methodology}\label{sec:alg}
In this section, we propose a 3-stage solution algorithm to solve the established RCVRPTW-TP-R, which includes the grouping approach, the route construction heuristic, and the scheduling approach. The method is inspired by the "route-then-schedule" scheme introduced by \cite{Luo_2022}. However, in \cite{Luo_2022}, each truck is assumed to serve only one pre-determined customer node, so its origin-destination pair is unique and pre-determined. By contrast, in our problem, a truck can serve multiple customer nodes; the sequence to serve these customers is yet to be determined; and the route of each truck cannot be decided until their customers' assignment and visiting sequence are given. Therefore, our problem is much more complicated and thus more challenging to solve than the one in \cite{Luo_2022}. Given this, 
our solution approach anticipates simplifying the problem by first determining the customers to serve for each truck, second routing each truck, and last assigning trucks with suitable schedules.

\subsection{Grouping approach}
In the first stage, we propose a Knapsack dynamic programming approach to divide customer nodes into several groups, such that the customer nodes in the same group will be served by the same truck as long as it has sufficient capacity. To fasten the grouping process, we first conduct a time-window check to identify the customers that are time-window feasible to avoid infeasible grouping. Particularly, the time-window feasibility is defined below:

\begin{definition}\label{proptw}
    Given two customer nodes $i$ and $j$, and that their earliest arrival times and latest departure times are $T^{EA}_i, T^{EA}_j$ and $T^{LD}_i, T^{LD}_j$, respectively, we claim that they are the least time-window feasible if the following inequality holds:
    $$
    \max(T^{LD}_j - T^{EA}_i, T^{LD}_i - T^{EA}_j) \geq t_{i,j}
    $$
    where $t_{i,j}$ stands for the shortest path travel time between customer node $i$ and $j$ under the road network.
\end{definition}

The intuition behind this definition is that if the shortest path travel time between two customer nodes is larger than the allowed time window range, then these two customer nodes definitely cannot be served by the same truck. By applying the time-window check defined in Definition \ref{proptw}, we can sift out some pairs of customer nodes that obviously cannot be allocated to the same truck, thus reducing the solution space.

The detailed process is given in Algorithm \ref{alg:Customerstw}.

\begin{breakablealgorithm}
\caption{Grouping Customers Based on Time-Window Feasibility}
\label{alg:Customerstw}
\begin{algorithmic}[1]
\Require{$\textit{Customers}$}
\Ensure{$\textit{Grouped customers based on time-window feasibility}$}
\State Initialize $\textit{GroupList} \leftarrow \emptyset$
\State Initialize $\textit{AssignedCustomers} \leftarrow \emptyset$
\While{$\textit{AssignedCustomers}$ size is less than the number of customers}
    \State Initialize $\textit{NewGroup} \leftarrow \emptyset$
    \State Add the first unassigned customer to $\textit{NewGroup}$
    \State Add this customer to $\textit{AssignedCustomers}$
    \Repeat
        \State $\textit{AddFlag} \leftarrow$ false
        \For{each unassigned customer $i$ in $\textit{Customers}$}
            \If{$i$ is time-window feasible with all customers in $\textit{NewGroup}$}
                \State Add $i$ to $\textit{NewGroup}$
                \State Add $i$ to $\textit{AssignedCustomers}$
                \State $\textit{AddFlag} \leftarrow$ true
                \State Break
            \EndIf
        \EndFor
    \Until{$\textit{AddFlag}$ is false}
    \State Add $\textit{NewGroup}$ to $\textit{GroupList}$
\EndWhile
\State \Return{$\textit{GroupList}$}
\end{algorithmic}
\end{breakablealgorithm}

Specifically, the process of Algorithm \ref{alg:Customerstw} can be explained as follows.
\begin{itemize}[label={[\arabic*]}]
    \item[{[1]}] Create an empty list, and put the first customer that has not been assigned to this list.
    \item[{[2]}] Iteratively add non-assigned customers to the list if the non-assigned customers are time-window feasible to all existing customers in the list.
    \item[{[3]}] Repeat step 2 until no more customers can be added
    \item[{[4]}] Repeat steps 1-3 until all customers are assigned to a specific list.
\end{itemize}

The time complexity of Algorithm \ref{alg:Customerstw} is $O(n^2)$, where $n$ denotes the number of customer nodes. By checking the time-window feasibility, the customers from different groups in the GroupList cannot be served by the same truck, and thus we can consider each group one by one. Therefore, the next step of grouping is to iteratively dispatch trucks to serve the customers of each group, and this process can be parallelly computed as the decision process of different groups will not interfere with each other. The grouping procedure is essentially a Knapsack problem that can be solved by the dynamic programming approach, so we elaborate on its corresponding Knapsack capacity, item weight, cost, and value function. Specifically, the value function is originally designed by us to better fit our problem structure. 

Knapsack capacity: each truck will be viewed as a Knapsack, so it corresponds to the truck capacity.

Item weight: we want to categorize customer nodes into trucks with fixed capacity, so it corresponds to the customer demand size.

Cost: the shortest travel time from each customer node to the depot.

Value function: Denote $F$ as the value function of the dynamic programming approach, $F(i,\hat{q})$ as the value function stage at customer node $i$ and capacity level $\hat{q}$, where $i$ is ordered from 1 to the last item in the group, and $0\leq \hat{q}\leq Q$, where $Q$ is the truck capacity.
If $i> 1$ and $q_i \leq \hat{q}$, then the value of retro-back stage is
$$
F(i,\hat{q}) = \max\left\{F(i-1,\hat{q}),F(i-1,\hat{q}-q_i)+(t^{max}+1-t_{i-1,i})\right\}
$$
where $t^{max}$ corresponds to the maximum value of all shortest paths of all pairs of nodes, and $t_{i,i-1}$ denotes the shortest path travel time between customer node $i$ and $i-1$, obtained by Floyd-Warshall algorithm with time complexity of $O(n^3)$. The retro-back updating equation illustrates that the function value of the current stage takes the greater value of two cases: the value of the previous stage with the same capacity when the current customer is not assigned or the value of the previous stage when the current customer is assigned plus the shortest path travel time between these two consecutive nodes. 
If $q_i > \hat{q}$, then it means that the current node cannot be added to the truck, i.e.,
$$
F(i,\hat{q}) = F(i-1,\hat{q})
$$

\noindent and we have the initial stages for the value function, when the capacity is 0 or when there is no customer not assigned.
$$
F(0,\hat{q}) = F(i,0) = 0
$$

Additionally, since we may dispatch multiple trucks to serve customer demands, the problem is equivalent to a multi-Knapsack problem. Therefore, we iteratively put customers into different trucks, and each truck will take the customers that maximize the value function while neglecting the influences of trucks considered later, so our approach is an essentially greedy approach. As a result, the order of customer nodes matters and will influence the solution quality, and thus we introduce perturbation and adjustment processes in Subsection \ref{sp} to weaken such an influence. Algorithm \ref{alg:dynamicprogramming_assignment} provides the pseudocode for solving the above dynamic programming problem with time complexity of $O(n^3Q)$.

\begin{breakablealgorithm}
\caption{Dynamic Programming for Grouping and Truck Dispatch}
\label{alg:dynamicprogramming_assignment}
\begin{algorithmic}[1]
\Require{$\textit{Dist}$, $\textit{Customers in GroupList}$, $\textit{Trucks}$}
\Ensure{$\textit{Customer assignment for each truck}$}
\State Initialize $\textit{CustomerAssignment}[v] \leftarrow \emptyset$ for each $v \in \textit{Trucks}$
\State Initialize $\textit{ServedCustomers} \leftarrow \emptyset$
\While{$\textit{ServedCustomers}$ size is less than the number of customers}
    \For{each $v \in \textit{Trucks}$}
        \State Initialize $F[i, \hat{q}] \leftarrow 0$ for all unserved customers $i$, $0 \leq \hat{q} \leq Q$.
        \For{each unserved customer $i$ in $\textit{Customers}$, ordered from $1$ to $n$}
            \For{$\hat{q} = 0$ to $Q$}
                \If{$q_i \leq \hat{q}$}
                    \State $F[i,\hat{q}] \leftarrow \max(F[i-1,\hat{q}], F[i-1, \hat{q}-q_i] + (t^{max} + 1 - t_{i-1,i}))$
                \Else
                    \State $F[i,\hat{q}] \leftarrow F[i-1, \hat{q}]$
                \EndIf
            \EndFor
        \EndFor
        \State $\textit{OptimalGroup} \leftarrow$ backtracking from $F[n, Q]$ to form a group of customers.
        \If{all customers in $\textit{OptimalGroup}$ can be served by truck $v$}
            \State $\textit{CustomerAssignment}[v] \leftarrow \textit{OptimalGroup}$
            \State Add customers in $\textit{OptimalGroup}$ to $\textit{ServedCustomers}$
        \EndIf
    \EndFor
\EndWhile
\State \Return{$\textit{CustomerAssignment}$}
\end{algorithmic}
\end{breakablealgorithm}
where the input $\textit{Dist}$ refers to an adjacency matrix that records the shortest path travel times between any two nodes of the network. Input $\textit{Trucks}$ are empty lists that customers will be allocated to, and each customer can be allocated to at most one truck. 


\subsection{Route construction heuristic}
In this subsection, we will construct routes for each dispatched truck by a modified insertion heuristic to minimize energy consumption individually. By intuition, each truck will tend to serve the customer node that is nearest to the depot and has the most demand first, and thus the truck can travel the remaining itinerary with a minimum weight, thereby saving more energy. However, in most cases, the customer node that is nearest to the depot may not have the most demand. Therefore, to take these two factors into consideration simultaneously, we propose a product, which is $(\eta Q+\gamma -\eta q_i)t_{0,i}$ for any customer node $i$, to estimate their importance. The customer node with a smaller product value will be considered more important, and it will be served first to reduce the total traveling cost. The modified insertion heuristic proposed in this paper stems from such an idea, and detailed procedures of the heuristic are provided in Algorithm \ref{alg:route_construction},. 
In the algorithm, $p_h$ represents the remaining load on the truck after departing from node $h$. 

\begin{breakablealgorithm}
\caption{Modified Insertion Heuristic for Route Construction}
\label{alg:route_construction}
\begin{algorithmic}[1]
\Require{$\textit{Customers}$ for each truck, $Q$, $q_i$, $t_{i,j}$, $\eta$, $\gamma$}
\Ensure{$\textit{Constructed routes}$ for each truck}
\For{each $\textit{Truck}$}
    \State Initialize $\textit{subTour} \leftarrow \emptyset$
    \State Select $i$ in current truck with smallest $(\eta Q+\gamma -\eta q_i)t_{0,i}$, generate $\textit{subTour} = (0-i-0)$
    \While{there are unvisited nodes in the current truck}
        \State Select $r$ not in $\textit{subTour}$ with smallest $(\eta Q+\gamma- \eta q_r - \eta \sum_{i\in subTour}q_i)t_{sub,r}$
        \State Find $(j,h)$ in $\textit{subTour}$ where $t_{j,r}(\eta p_h+\eta q_r+\gamma) + t_{r,h}(\eta p_h+\gamma) - t_{j,h}(\eta p_{h}+\gamma)$ is smallest
        \State Insert $r$ into $\textit{subTour}$ between $j$ and $h$
    \EndWhile
    \State Add $\textit{subTour}$ to the routes of the current $\textit{Truck}$
\EndFor
\State \Return{$\textit{Constructed routes}$ for each truck}
\end{algorithmic}
\end{breakablealgorithm}

The algorithm's time complexity is $O(n^3)$, and it can be explained as follows: 
\begin{enumerate}[label={[\arabic*]}]
    \item Pick a customer node $i$ that will be served by the current truck and with the smallest value of $(\eta Q+\gamma -\eta q_i)t_{0,i}$, to generate a subtour $(0-i-0)$.
    \item Pick any customer node $r$ not in the current subtour and with the smallest value of $(\eta Q+\gamma- \eta q_r - \eta \sum_{i\in sub}q_i)t_{sub,r}$, where $\sum_{i\in sub}q_i$ indicates the total customer demand served by the current subtour, and $t_{sub,r}$ represents the minimum value of the shortest travel time from node $r$ to any customer node or depot in the current subtour.
    \item Insert the node $r$ into the current subtour, find a road segment $(j,h)$ such that $t_{j,r}(\eta p_h+\eta q_r+\gamma) + t_{r,h}(\eta p_h+\gamma) - t_{j,h}(\eta p_{h}+\gamma)$ is the smallest. Since $p_h$ represents the remaining load after leaving node $h$, the whole expression basically calculates the increase in energy consumption in the road segment between node $j$ and $h$ when customer node $r$ is inserted between them. To be mentioned, $p_h$ can be easily calculated based on the current subtour.
    \item Iteratively repeat steps 2 and 3 until all customer nodes of the current truck are added to the subtour.
\end{enumerate}

To be mentioned, $\eta$, the marginal fuel consumption coefficient, is used in our calculation to precisely capture the influences of weight on fuel consumption. In addition, in the first step, we pick $(\eta Q+\gamma -\eta q_i)t_{0,i}$ to consider the weight and distance simultaneously. Supposing that the truck departs from the depot with a full load, taking the smallest value of $(\eta Q+\gamma -\eta q_i)t_{0,i}$ will encourage the truck to first serve the customer node $i$ that both has a small distance to the depot, $t_{0,i}$, and large demand $\eta q_i$. Accordingly, our heuristic will most likely put the most influential customer node in the first place of the subtour. In the following steps, we use a similar idea to iteratively add each customer node into the subtour and finally construct a complete round trip for each dispatched truck.

\subsection{Scheduling approach} \label{sp}
As we mentioned above, we omit time-window constraints to maximize platoon formations and minimize truck detours by allowing nearby customer nodes in the same group. The solution obtained from the first two stages can act as an upper bound to the RCVRPTW-TP-R problem if it is feasible for the scheduling problem, and it maximizes savings and benefits. However, some customer nodes may be time-infeasible given the truck delivery routes provided by the previous two stages. Therefore, we need to correct them in this stage and output a feasible solution for the overall model.

The way to correct time-infeasible customer nodes is based on the ``route-then-schedule" concept from \cite{Luo_2022}. Our group approach and route construction heuristic can be seen as a way to solve the routing problem in the ``route-then-schedule" loop. For the scheduling part, we construct a scheduling problem that maximizes platoon fuel savings. The scheduling problem is a pruned version of the original problem as it only investigates the constructed routes instead of the whole network, and the customers to be served for each truck are pre-determined and cannot be changed. Therefore, the solution space of the scheduling problem is restricted to a small-scale problem and can be solved fast. The detailed formulation of the scheduling problem is provided as follows:

\begin{flalign}
&\min  \sum_{(i,j)\in \cup_k \mathcal{R}_k} \dfrac{\alpha}{\gamma} t_{i,j}\left(\sum_{k\in \mathcal{K}_{i,j}} \left(v^{(l)}_{i,j,k} + (1-\beta) v^{(f)}_{i,j,k}\right)\right) \notag\\
&\mbox{s.t.} \notag \\
&\sum_{k_2\in  \mathcal{K}_{i,j}:k_2 \neq k} f_{i, j, k, k_2} \leq (L - 1) \times l_{i, j, k}
&\forall (i,j)\in \mathcal{R}_{k}, k\in K,&\notag\\
&\sum_{k_2\in  \mathcal{K}_{i,j}:k_2 \neq k} f_{i, j, k_2,k} + l_{i, j, k} \leq 1
&\forall (i,j)\in \mathcal{R}_{k}, k\in K, &\notag\\
&-M(1 - f_{i, j, k_1, k_2}) \leq s_{i, k_1} + w_{i, k_1} - s_{i, k_2} - w_{i, k_2}
&\forall (i,j)\in \cup_{k\in K}\mathcal{R}_k,  k_1\neq k_2\in \mathcal{K}_{i,j}, &\notag\\
&s_{i, k_1} + w_{i, k_1} - s_{i, k_2} - w_{i, k_2} \leq M(1 - f_{i, j, k_1, k_2})
&\forall (i,j)\in \cup_{k\in K}\mathcal{R}_k,  k_1\neq k_2\in \mathcal{K}_{i,j}, &\notag\\
&s_{i, k} + t_{i, j} + w_{i, k} - M(1 - x_{i, j, k}) \leq s_{j, k}
&\forall (i,j)\in \mathcal{R}_k,  k\in K, &\notag\\
&T^{EA}_i \times g_{i, k} \leq s_{i, k}
&\forall i\in C, k\in K, &\notag\\
&s_{i, k} + w_{i, k} \leq T^{LD}_i \times g_{i, k}
&\forall i\in C, k\in K, &\notag\\
&v^{(l)}_{i,j,k} \leq (\eta Q+\gamma)l_{i,j,k}
&\forall (i,j)\in \cup_{k\in K}\mathcal{R}_k,k\in K, &\notag\\
&v^{(l)}_{i,j,k} \leq \eta y_{i,k} + \gamma 
&\forall (i,j)\in \cup_{k\in K}\mathcal{R}_k,k\in K, &\notag\\
&v^{(l)}_{i,j,k} \geq \eta y_{i,k} + \gamma - (\eta Q+\gamma)(1-l_{i,j,k})
&\forall (i,j)\in \cup_{k\in K}\mathcal{R}_k,k\in K, &\notag\\
&v^{(f)}_{i,j,k} \leq (\eta Q+\gamma)\sum_{k_1:k_1\neq k}f_{i,j,k_1,k}
&\forall (i,j)\in \cup_{k\in K}\mathcal{R}_k,k\in K, &\notag\\
&v^{(f)}_{i,j,k} \leq \eta y_{i,k} + \gamma 
&\forall (i,j)\in \cup_{k\in K}\mathcal{R}_k,k\in K, &\notag\\
&v^{(f)}_{i,j,k} \geq \eta y_{i,k} + \gamma - (\eta Q+\gamma)(1-\sum_{k_1:k_1\neq k}f_{i,j,k_1,k})
&\forall (i,j)\in \cup_{k\in K}\mathcal{R}_k,k\in K, &\notag\\
&l_{i, j,k}\in \left\{0,1\right\}
&\forall i,j\in N,k\in K, &\notag\\
&f_{i,j,k_1,k_2}\in \left\{0,1\right\}
&\forall (i,j)\in \cup_{k\in K}\mathcal{R}_k,k_1,k_2\in K,k_1\ne k_2, &\notag\\
&s_{i, k}\geq 0
&\forall i\in N, k\in K, &\notag\\
&w_{i, k}\geq 0
&\forall i\in N, k\in K, &\notag\\
&v^{(l)}_{i,j,k} \geq 0
&\forall (i,j)\in \cup_{k\in K}\mathcal{R}_k, k\in K, &\notag\\
&v^{(f)}_{i,j,k} \geq 0
&\forall (i,j)\in \cup_{k\in K}\mathcal{R}_k, k\in K, &\notag
\end{flalign}

Particularly, the decision variables defined in the original problem, i.e., RCVRPTW-TP-R, including $x_{i,j,k}$, $y_{i,k}$ and $g_{i,k}$, will not be viewed as variables in this sub-problem. Instead, they will be treated as parameters, and their values can be determined  ($x_{i,j,k} = 1$ if $(i,j)\in \mathcal{R}_{k}$, and $\mathcal{K}_{i,j}$ collects all the trucks that travel on link $(i,j)$) based on the route assignment of trucks, which are already decided in the previous two stages, before solving the scheduling problem in each iteration. It can be further witnessed that the scheduling problem uses a partial objective function, which is the second term of the objective function in RCVRPTW-TP-R, as the number of trucks to be dispatched has already been decided and we only care about maximizing the platooning benefit. The other constraints that need to be considered are platoon constraints (\ref{8})-(\ref{12}), time window constraints (\ref{13}) and (\ref{14}), integrity constraints (\ref{15})-(\ref{21}), and constraints (\ref{22})-(\ref{29}) that used to linearize nonlinear terms. In addition, constraints considered in the scheduling problem only consider edges $(i,j)$ that are covered by the route assignment obtained in the previous two stages, denoted as the union of all traveled routes, $\cup_{k\in K} \mathcal{R}_k$. 

However, the scheduling problem may appear to be infeasible because we have omitted time constraints, and the least time-window infeasibility check declared in Definition \ref{proptw}, though eliminating many obvious infeasible grouping selections, may be insufficient in the instances with rather strict time windows for some customers. Therefore, in this case, we will return to the grouping stage and shuffle the order of all customer nodes, and then re-generate the truck assignments and reconstruct the route of each truck.

\begin{figure}[htbp]
	\centering
    \includegraphics[scale=0.9]{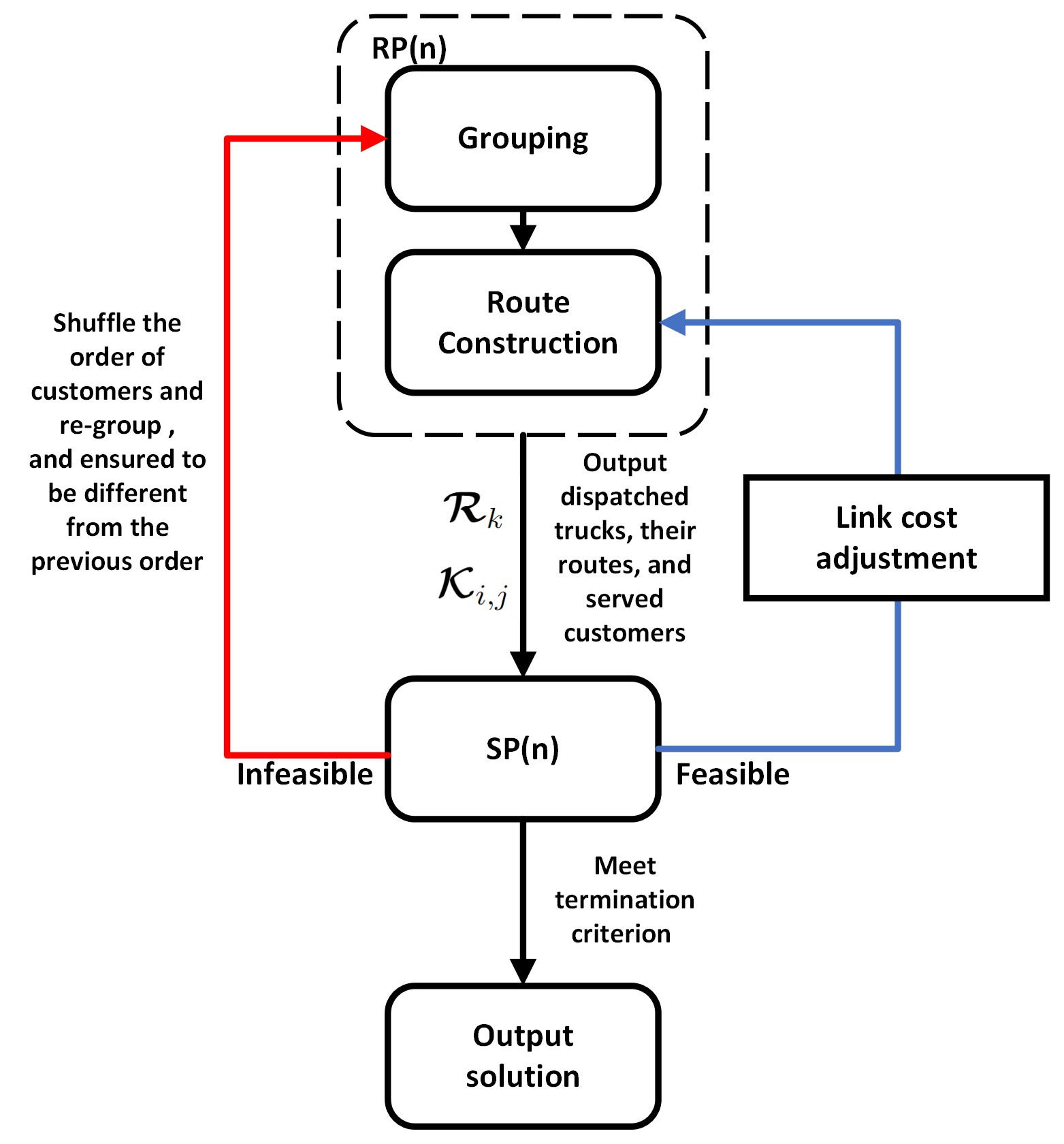}
	\caption{Iterative process of algorithm}
    \label{fig:alg}
\end{figure}

On the other hand, when the scheduling problem is feasible, then the optimal solution obtained from the scheduling problem acts as a feasible solution to the original problem and also provides an upper bound. In this case, adjustment procedures will be applied to improve the current solution. Therefore, the overall procedures of our iterative solution algorithm can be represented as a diagram that is shown in Figure \ref{fig:alg}.

To be mentioned, we denote the routing problem at iteration $n$ as $RP(n)$, and the scheduling problem at iteration $n$ as $SP(n)$, and the termination criterion is met when the time limit occurs or the current truck routes' of the solution has not changed for many iterations. A detailed explanation of link cost adjustment is provided below.

\vskip 0.2in

\noindent \textbf{Modified link cost}\par

In Figure \ref{fig:alg}, we modify the link costs in the feedback loop when the scheduling problem is feasible, and the modified link costs will act as an input to the route construction stage of the routing problem in the next iteration. Accordingly, the route construction stage will generate routes for each truck based on the network with modified link costs.  Therefore, we define the link cost of truck $k$ that inputs to the routing problem at iteration $n$ as $t_{i,j,k}^{(n)}, \forall (i,j)\in A,k\in K$, and define the traversed routes of trucks at iteration $n$ as $\cup_{k\in K} \mathcal{R}_k^{(n)}$. 
Assuming that the $SP(n)$ is feasible at iteration $n$, then the presumed link traveling cost (in time) for assigning a truck $k_1\in K$ to an edge $(i,j)\in A$ at iteration $(n+1)$ during the route construction stage will be updated based on the following schemes \citep{Luo_2022}:
\begin{equation}\label{update}
t_{i,j,k_1}^{(n+1)}=\left\{\begin{array}{ll}
\frac{t_{i, j}^{\text {platoon }}\left(\left|\mathcal{P}_{i,j,k_1}^{(n)}\right|\right)}{\left|\mathcal{P}_{i,j,k_1}^{(n)}\right|} & \text { if }(i, j) \in \mathcal{R}_{k_1}^{(n)}, \tag{30}\\
\left(1-\beta\right)t_{i, j} & \text { if }(i, j) \in  \cup_{k \in K} \mathcal{R}_k^{(n)} \backslash \mathcal{R}_{k_1}^{(n)} , \\
t_{i,j} & \text {otherwise},
\end{array}\right.
\end{equation}
where $\left|\mathcal{P}_{i,j,k_1}^{(n)}\right|$ denotes the number of trucks in the same platoon with truck $k_1$ on edge $(i,j)$ at iteration $n$, so it equals to 1 if truck $k_1$ is traveling alone. The term $t_{i, j}^{\text {platoon }}\left(\left|\mathcal{P}_{i,j,k_1}^{(n)}\right|\right)$ represents the traveling cost of the whole platoon, while the load of each truck is assumed to be 0 during calculation.

The intuition of the updating criterion in (\ref{update}) can be explained as follows. In the first case, the truck is traveling the link that was traversed by it in the previous iteration, so we assume that it will incur the platoon-averaged traveling cost. In the second case, when the truck is traveling on an edge that was not explored by it in the previous iteration, but this edge was traversed by another truck before, then the adjusted link cost $\left(1-\beta\right)t_{i, j}$ aims to reflect that the truck can travel as a following truck on this edge. The third case shows that if the truck traverses the edge that was not explored by any truck, then we directly use the original link traveling cost. By such a link cost modification scheme, we can encourage trucks to explore different routes and form a platoon in the next iteration.

\section{Numerical Experiments}\label{sec:exp}
In this section, experiments based on a virtual small network and a real-world road network, the Yangtze River Delta network, will be conducted to examine the performance of our proposed model and solution algorithm. All experiments were programmed by Java with Intellij interface and were executed on a PC with 6-core 12-thread Intel(R) i7-8700k 3.7GHz CPU and 32 GB 3000MHz RAM. To be mentioned, the termination criterion is met when the time limit (3600s) occurs or the current truck routes' of the solution has not changed for many iterations (we chose 10).

The default parameter settings are provided in Table \ref{tab:par}. According to the real-world data \citep{salary}, the average base salary for truck drivers in the United States per day is \$271, so we set $c_1 = 271$. Additionally, the energy consumption rate, $\alpha$, given that electric trucks are not prevalent yet, we adopt the fuel consumption rate of trucks for current experiments. 
As indicated by \cite{Fuelconsump}, the fuel consumption rate of a truck is 33 liters per 100 kilometers for a Class 8 truck, and the average traveling speed of trucks in the U.S. on highways is 55 mph, so the fuel consumption rate is calculated as 29.2 liters per hour. Therefore, the fuel consumption rate is calculated as $29.2\times 1.05 = 30.7$ dollar per hour, based on the average gasoline price per liter in the U.S. in 2022. Furthermore, the weight of a Class 8 truck's tractor ranges from 8 tons to 10 tons, and the maximum weight limit of a Class 8 truck on U.S. highways is 80,000 pounds, which is equivalent to 36.2 tons. On the other hand, the maximum weight limit on China highways is 31 tons including truckload. Therefore, to simplify the calculation, we take 10 tons as the static weight of trucks, $\gamma$, and the maximum capacity, $Q$, to be 20 so that the maximum possible gross vehicle weight will not exceed 30 tons. In addition, according to reports and studies on truck fuel consumptions \citep{weighteffect, weighteconomy, quantitativeeffect}, the fuel consumption rate at a full-load state, which equals the truck capacity $Q$, can be 15\%-25\% larger than the basic rate, so we choose 20\% in our experiments. Accordingly, it is equivalent to saying that the term $\dfrac{\eta y +\gamma}{\gamma}$ should be 1.2 when the load is full, i.e., $y = Q = 20$. Therefore, the marginal fuel consumption coefficient, $\eta$, can be calculated as $\dfrac{1.2\times \gamma - \gamma}{y} = \dfrac{1.2\times 10-10}{20} = 0.1$ per every additional unit of load (unit in ton). Lastly, as \cite{drayage2020} suggested, the platoon saving ratio of the following vehicles in a platoon may vary from 0.05 to 0.15, so we set $\beta = 0.1$ for the base cases, and the maximum platoon size, $L$, is selected to be 4 to enable higher platoon potential. 
\begin{table}[htbp]
  \centering
  \caption{Parameter Setting}
  \label{tab:par}
  \begin{tabular}{lll}
    \toprule
     Parameter & Description & Value\\
    \hline
    $c_1$ & Vehicles' dispatch cost per truck & \$271\\
    $c_2$ & Weight coefficient to balance two kinds of costs & 1\\
    $\alpha$ & Base fuel consumption rate respect to unit time (truck static weight) & \$30.7/hr\\
    $\gamma$ & Static weight of trucks & 10 (tons)\\
    $Q$ & Maximum truck capacity & 20 (tons) \\
    $\beta$ & Platoon's fuel cost savings ratio & 0.1\\
    $L$ & Maximum platoon size  & 4\\
    $\eta$ & Marginal fuel consumption coefficient for a unit of additional weight & 0.1\\
    \bottomrule
\end{tabular}
\end{table}

\subsection{A small scale experiment}

In this subsection, we manually design a small test case with 8 customer demands on a 4 by 4 grid network. The network size is refined to ensure that an exact solution can be obtained by commercial solvers, such as CPLEX, in a timely manner. The test case is solved by our provided model when the platooning is enabled and also solved when the platooning is prohibited, reflecting the traditional delivery service. The latter scenario can be simply accomplished by our model when the platoon size limit is set to 1. Results comparison between these two scenarios will be presented in the remainder of this subsection to demonstrate the benefits of utilizing truck platooning technology in a capacitated delivery service.

\begin{figure}[htbp]
	\centering
    \includegraphics[scale=0.7]{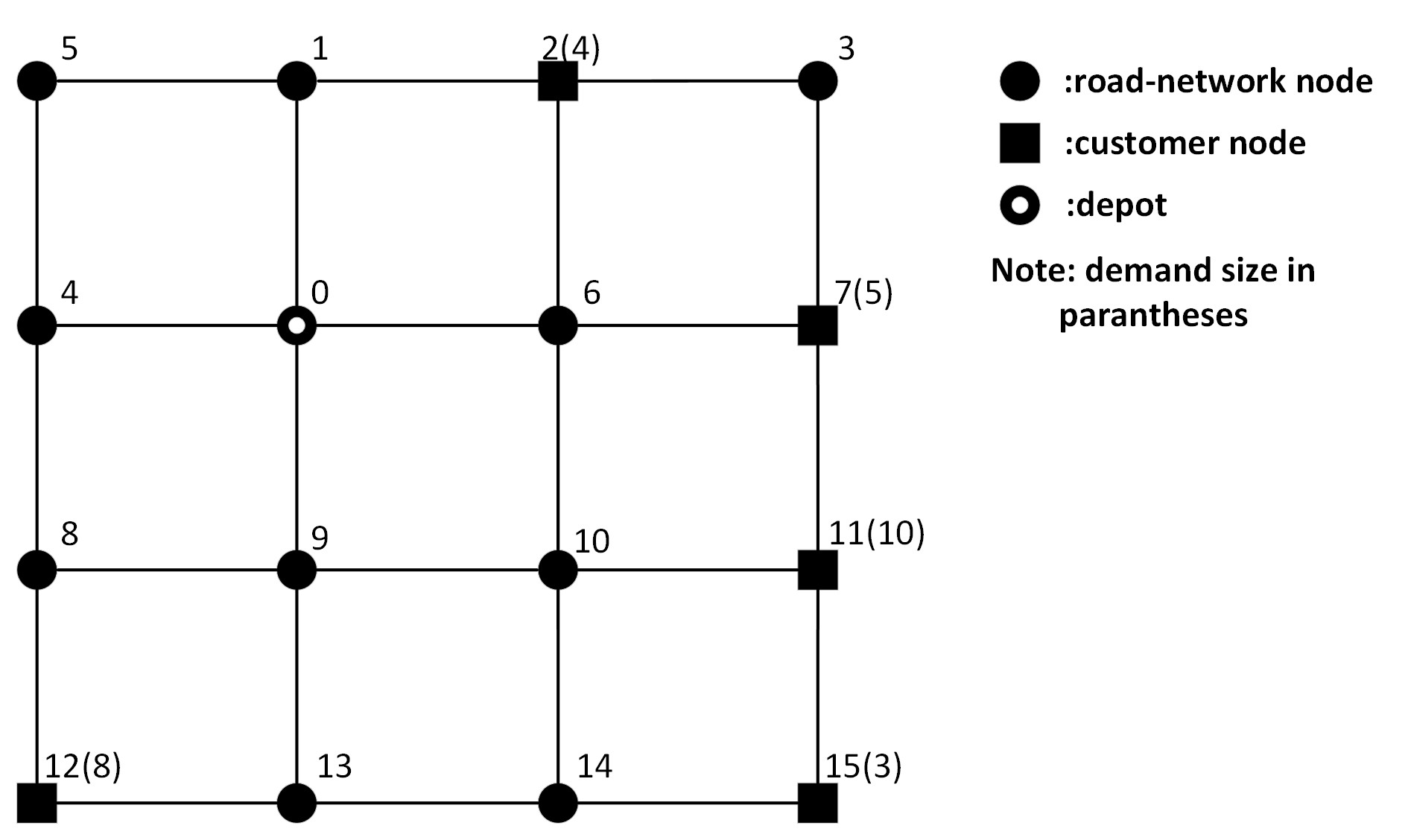}
	\caption{A small scale network}
    \label{fig::grid}
\end{figure}

The test network is shown in Figure \ref{fig::grid}, and it is a 4 by 4 grid network with identical edges. Each edge has a length of 3, representing that it takes 3 hours for a truck to travel from one node to an adjacent node. We separate different kinds of nodes by distinct shapes: square nodes are customer nodes; solid round nodes are road network nodes, and the only circle node represents the central depot. The label of each node is located at the top right corner of each one, such as node 0 as the central depot. For each customer node, the size of the customer demand (in tons) is written within the brackets. To maximize the platoon potential and operation difference between the two scenarios, we make the service time window of each customer node large enough ([0,36]) so that it can improve the solution results. 

Table \ref{operacomp} illustrates the operation cost for both scenarios. We can witness that, when platooning is prohibited, referring to the traditional delivery service, the total operation cost is a bit higher than the scenario when platooning is allowed. In this test case, the dispatch cost for both scenarios is the same, indicating that the numbers of dispatched drivers are the same for both scenarios. Therefore, the difference in operation cost all comes from the energy savings by platooning. The total savings in percentage can be therefore calculated: $\dfrac{1375.8-1362.3}{1375.8} = 0.98\%$. The amount of savings is not quite significant in this case because platoons with only two trucks are formed, and customer nodes are scattered so that platooning links are limited. In real instances, larger platoons can be formed, and there exist more common paths among different truck routes, resulting in more significant platoon savings. Furthermore, even such a marginal percentage increase can translate into substantial cost savings for the trucking industry, considering its considerable operational expenses in practice.

\begin{figure}[htbp]
	\centering
	\subfloat[Maximum platoon size is 4]{\label{fig:gridplatoon}\includegraphics[scale=0.7]{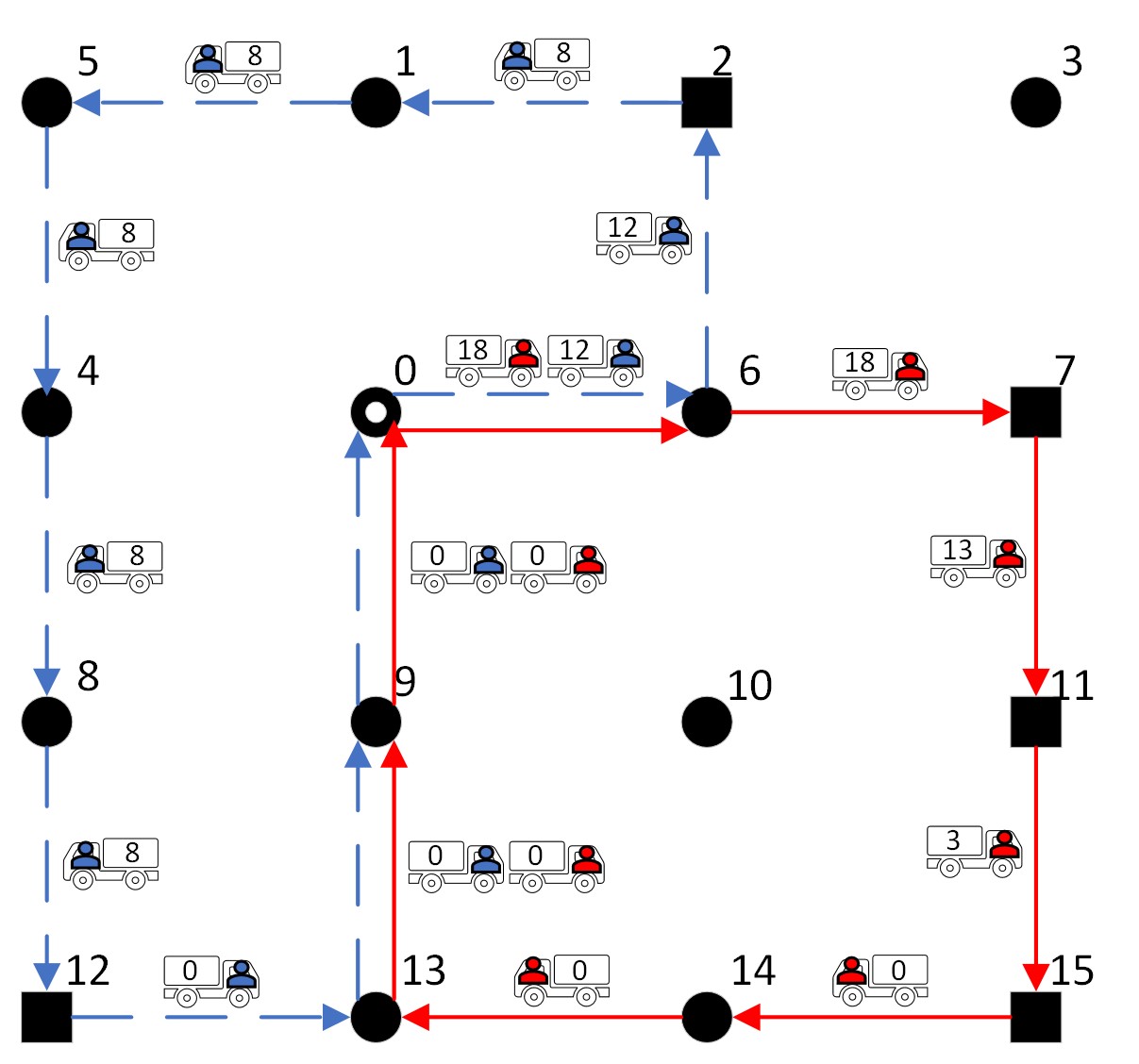}}\quad
    \subfloat[Platooning is prohibited]{\label{fig:gridvrp}\includegraphics[scale=0.7]{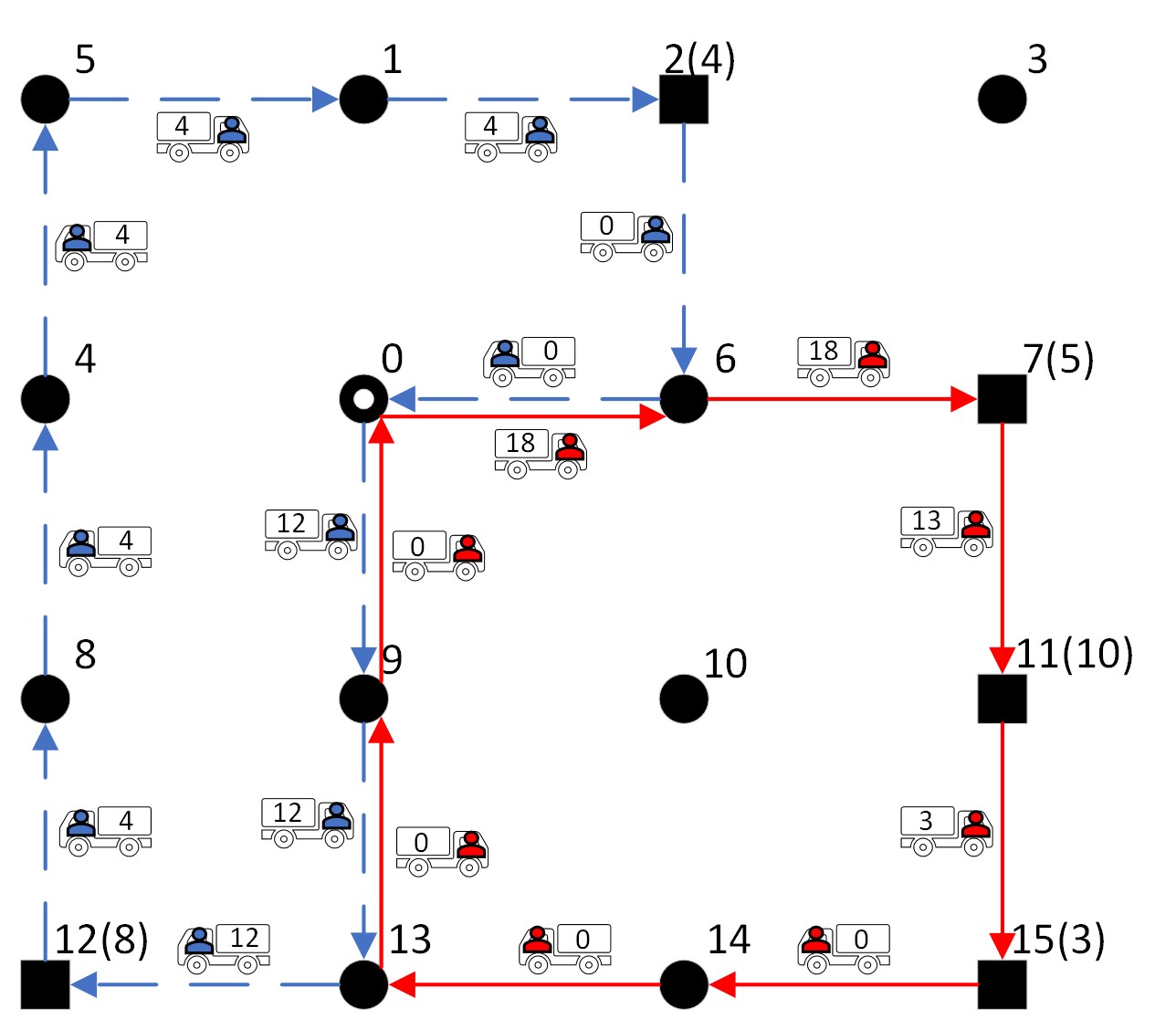}}\\
	\caption{Operation plans with and without platooning}
    \label{fig::gridresult}
\end{figure}

We further discuss the impact of platooning on operation plans, and visually present the operation plans under two scenarios in Figure \ref{fig::gridresult}. When platooning is not allowed, then the original problem is equivalent to a multiple-vehicle routing problem since the order to serve customers and the route to serve customers for each truck will no longer be influenced by other trucks' routes. The resulting operation plan under this scenario is provided in Figure \ref{fig:gridvrp}, and the routes of different trucks are drawn with different colors and line formats. We further record the load of each truck along its traveling route. It can be seen that two trucks are required to fulfill the total demand. On the other hand, the operation plan when platooning is allowed is shown in Figure \ref{fig:gridplatoon}. Two trucks are also required under this scenario, and the truck that takes the red solid line has exactly the same route as the other scenario. However, the other truck is traveling on an order: 0-6-2-1-5-4-8-12-13-9-0 under the platooning scenario, which is exactly reversed as 0-9-13-12-8-4-5-1-2-6-0 under the no-platooning scenario, in order to form a platoon on links (0,6), (13,9), and (9,0). Additionally, the truck with a higher load is acting as the following truck on the link (0,6) to maximize platoon savings, reflecting that the truck weight does influence the platooning formations. Furthermore, the truck that travels on the blue dashed line visits customer node 12 ahead of customer node 2 when platooning is prohibited because the delivery of heavier loads first can induce more energy savings. However, this truck visits customer node 2 first under the platooning case in order to form a platoon with the other truck, saying that platooning does influence the order to serve customer nodes and the routes adopted by trucks.

\begin{table}[htbp]
\caption{Comparison of operation cost between the without platooning scenario and platooning scenario}
\begin{tabular}{lcc}
\toprule

\textbf{Cost type} &  Without platooning & With platooning\\ 
\hline
\midrule
\textbf{Total operation cost} & 1375.8 & 1362.3  \\
\textbf{Dispatch cost} & 542 & 542\\
\textbf{Energy consumption cost} & 833.8 & 820.3 \\
\hline

\bottomrule
\end{tabular}
\label{operacomp}
\end{table}

\subsection{Experiments on the Yangtze River Delta network}
\begin{figure}[htbp]
	\centering
    \includegraphics[scale=0.12]{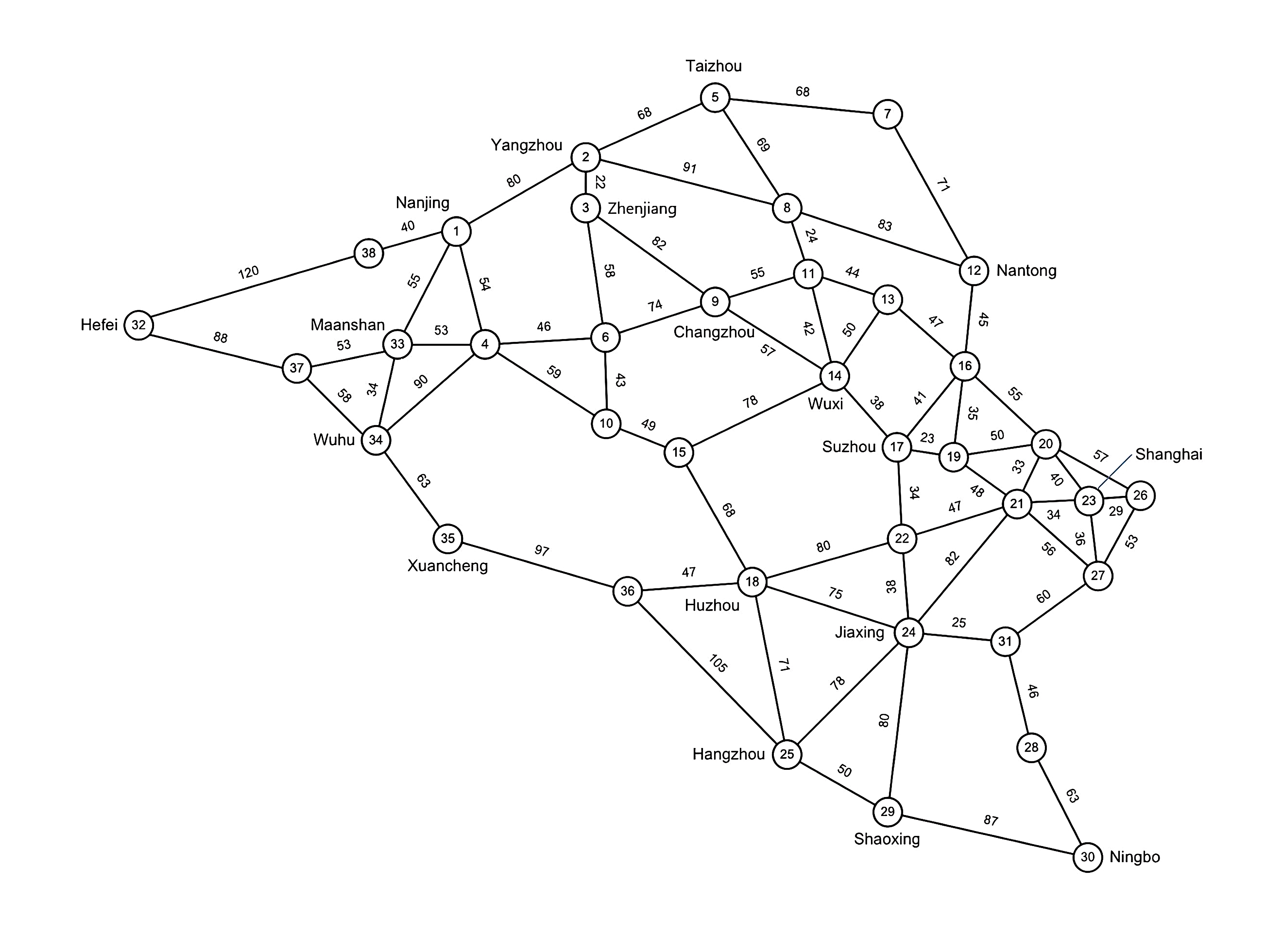}
	\caption{Yangtze River Delta network \citep{delta2022}}
    \label{fig::yang}
\end{figure}

In this subsection, we provide a numerical experiment on a real-world transportation network, the Yangtze River Delta network. It is located in the central east area of China, an economic center where logistics play an essential role. The network consists of 38 nodes or cities, and the whole network is presented in Figure \ref{fig::yang} with distances between nodes labeled. The units of labeled distances are in kilometers, so we manually transform them into time by dividing an average truck speed on highways, which is 88.5km per hour as also mentioned before. 

For the test case design, we set node 1, Nanjing, as our depot, and randomly chose some other nodes in the network as customer nodes with their customer demands randomly assigned, varying from 1 to 10 tons. Time windows are also randomly generated, and we set the latest departure times of each customer to be at least 20 (representing 20 hours) larger than its earliest arrival time to enable more platoon potential. In addition, we also ensured that the latest departure time of each customer node is larger than or equal to the shortest travel time between the depot and the customer node to guarantee the problem's feasibility.  We tested our results on test cases of various sizes, each with a different number of customer nodes, ranging from 5 to 20. Each test case will be executed by both CPLEX and our proposed algorithm, and we set a large time limit, which is 10,800 seconds, for CPLEX to output a high-quality solution while not wasting too much time. The running time and result comparisons between CPLEX and our proposed algorithm are presented in Table \ref{Yangresult}. In addition, for relatively small test cases, with the number of customer nodes ranging from 5 to 10, 5 instances were randomly generated for each case to provide a more thorough and even comparison between CPLEX and heuristic's performances.

\begin{table}[htbp]
\caption{Solution results}
\label{Yangresult}
\centering
\begin{threeparttable}
\resizebox{\columnwidth}{!}{\begin{tabular}{llllllll}
\toprule
\multicolumn{2}{l}{\textbf{Instance}}& \multicolumn{6}{l}{\textbf{Result}} \\
\cmidrule(lr){1-2} \cmidrule(lr){3-8}
\makecell{\textbf{Number of} \\ \textbf{customer nodes}} & \textbf{Execution count}  & \makecell{\textbf{CPLEX solution}} & \makecell{\textbf{CPLEX's} \\ \textbf{CPU time}}  & \textbf{Heuristic solution} & \makecell{\textbf{Heuristics's} \\ \textbf{CPU time}}   &\textbf{Gap} &\textbf{platooning benefit}\\ 
\hline
\midrule
5 & 1 &1,271 & 10,302s & 1,271 & 4.1s & 0 & 3.5\\
& 2& 1,672& 9,610s & 1,672& 7.3s & 0 & 6.6\\
& 3& 1,678& 9,885s & 1,678 & 5.6s & 0 & 0\\
& 4& 1,461& 10,152s & 1,461 & 6.9s & 0 & 0\\
& 5& 1,534& 10,800s & 1,534 & 6.8s & 0 & 8.0\\
6 & 1& 1,575 & 10,800s & 1,582 & 10.4s & 0.44\% & 0\\
& 2& 1,823& 10,800s & 1,819 & 12.7s & -0.22\% & 6.8 \\
& 3& 1,826 & 10,800s & 1,839 & 7.4s & 0.71\% & 7.7\\
& 4& 1,672& 10,800s & 1,672 & 8.4s & 0 & 0\\
& 5& 1,845& 10,800s & 1,845 & 7.6s & 0 & 0\\
7 & 1&1,928 & 10,800s & 1,882 & 23.2s & -2.39\% & 10.1\\
& 2& 2,025& 10,800s & 2,011& 22.7s & -0.69\% & 8.6\\
& 3& 1,839& 10,800s & 1,833 & 19.8s & -0.33\% & 6.9\\
& 4& 2,244& 10,800s & 2,098 & 25.3s & -6.5\% & 9.9\\
& 5& 2,541& 10,800s & 2,369 & 18.6s & -6.77\% & 0\\
8 & 1&3,234 & 10,800s & 3,055 & 30.7s & -2.87\% & 0\\
& 2& 2,539& 10,800s & 2,535 & 33.4s & -0.16\% & 21.5\\
& 3& 3,247& 10,800s & 3,085 & 56.1s & -4.99\% & 8.6\\
& 4& 2,560 & 10,800s & 2,535 & 47.3s & -0.98\% & 0\\
& 5& 2,828 & 10,800s & 2,693 & 57.5s & -4.77\% & 7.1\\
9 & 1&3,532 & 10,800s & 3,403 & 39.4s & -3.65\% & 13.5\\
& 2& 2,402& 10,800s & 2,249 & 43.9s & -6.37\% & 0\\
& 3& 3,821& 10,800s & 3,522 & 73.8s & -7.83\% & 23.3\\
& 4& N/A & 10,800s & 2,682 & 59.5s & N/A & 0\\
& 5& N/A & 10,800s & 2,838 & 62.1s & N/A & 33\\
10 & 1&N/A & 10,800s & 2,654 & 68.1s & N/A & 11.7\\
& 2& N/A& 10,800s & 2,845 & 61.1s & N/A & 18.8\\
& 3& 2,960& 10,800s & 2,858 & 55.2s & -3.45\% & 14.7\\
& 4& N/A& 10,800s & 3,104 & 70.5s & N/A & 6.8\\
& 5& 3,145& 10,800s & 2,990 & 53.4s & -4.93\% & 0\\
15 & 1&N/A & 10,800s & 4,128 & 112.3s & N/A & 34.7\\
20 & 1&N/A & 10,800s & 5,196 & 142.1s & N/A & 86.5\\
25 & 1&N/A & 10,800s & 6,083 & 154s & N/A & 113.9\\
\hline

\bottomrule
\end{tabular}}
\begin{tablenotes}
    \footnotesize
    \item Gap: the relative gap between CPLEX and heuristic results, and calculated by $\dfrac{Obj_{heuristic}-Obj_{CPLEX}}{Obj_{CPLEX}}$.
\end{tablenotes}
\end{threeparttable}
\end{table}

As we can see in Table \ref{Yangresult}, CPLEX always outputs the optimal solution in the cases in which the number of customer nodes equals 5, but an optimal solution cannot be guaranteed when the problem size grows larger. For instance, in the cases in which the number of customer nodes equals 6, our algorithm outperforms CPLEX by 0.22\%. In addition, when the problem size becomes large enough, such as the cases in which the number of customer nodes equals 9 and 10, CPLEX even fails to find a feasible solution within the time limit. Therefore, we did not perform multiple executions with CPLEX for even larger cases, when the number of customer nodes was 15, 20, or 25, because CPLEX could hardly solve the problem in a timely manner. Conversely, our proposed heuristic algorithm can yield solutions for all cases within 200 seconds, though an increase in solving time as the number of customer nodes scales can be witnessed. 
Although the number of all possible sequence orders of customers increases factorially as the number of customers scales, our solution algorithm manages the growth in solving time partially due to the time-window feasibility check. On the other hand, our proposed solution algorithm performs extremely well in small instances with 5 and 6 customer nodes, resulting in 0 and below 1\% gaps with the CPLEX optimal solution, respectively. Furthermore, it outperforms CPLEX for the cases with 7 and 8 customers with more than 5\% in several executions. Generally speaking, the benefits of our algorithm over CPLEX in producing a high-quality solution become increasingly significant as the problem becomes more complicated. 

We also provide a glance at the platooning benefit under each case. The platooning benefit is measured by the operation cost difference between the platooning case and the case without platooning. Therefore, we obtain the platooning benefit by running each case twice, once with the platoon size limit of 4, and once with the platoon size limit of 1 to represent the no-platoon case, and then calculate the cost differences between the two cases. We can expect that the main difference comes from the energy savings from the trailing trucks, and the other part due to improved dispatching strategy enhanced by truck platooning technology. In addition, due to the randomness of generated customer nodes, trucks may not platoon together in their optimal trajectories, and thus we can see 0 platooning benefit in some cases. Though the platooning benefit varies from case to case, we can still witness a rise in the platooning benefit as the problem size scales. For instance, the platooning benefit does not exceed 10 for the case of 5 customer nodes, and reaches a peak of 21.5 and 33 when the number of customer nodes gets to 8 and 9, respectively. The increase in platooning benefit is more notable when the problem size exceeds 15. 

\subsection{Sensitivity Analysis}
In this subsection, sensitivity analyses are conducted to quantitatively represent the impact of key parameters, including the fuel consumption rate, vehicle dispatch cost, and platoon size limit, on the results. All the experiments are designed on the Yangtze River Delta network, and the results are obtained from our proposed algorithm because CPLEX fails to provide even a feasible solution when the number of customer nodes exceeds 10.

In addition, to set up the experiments for sensitivity analyses, we constructed 5 different test scenarios, each with different customer demands and time windows, and problem feasibility always held during scenario design. Furthermore, each scenario has 4 cases, each with a different number of customer nodes, i.e., 5, 10, 15, and 20. We also ensured that the case with more customer nodes must include the customer nodes of the smaller case. For instance, the 20-customer case must include the 15 customers in the 15-customer case, and the 15-customer case must include the 10 customers in the 10-customer case, and so on. By doing so, the platooning benefit is comparable across cases of different sizes because the larger case always contains the customer demands of the smaller case. Accordingly, it became possible for us to investigate if platooning benefit is correlated with problem sizes. In addition, test cases will be identical in terms of demand information for each experiment by fixing the random seed. To sum up, we calculated the average results across the 5 test scenarios for each case with specific parameter settings and problem sizes to reduce the influence of randomness on our solution results.

\subsubsection{Fuel consumption rate}
\begin{figure}[htbp]
	\centering
	\subfloat[Operation cost]{\label{fig:fuela}\includegraphics[width=0.48\textwidth]{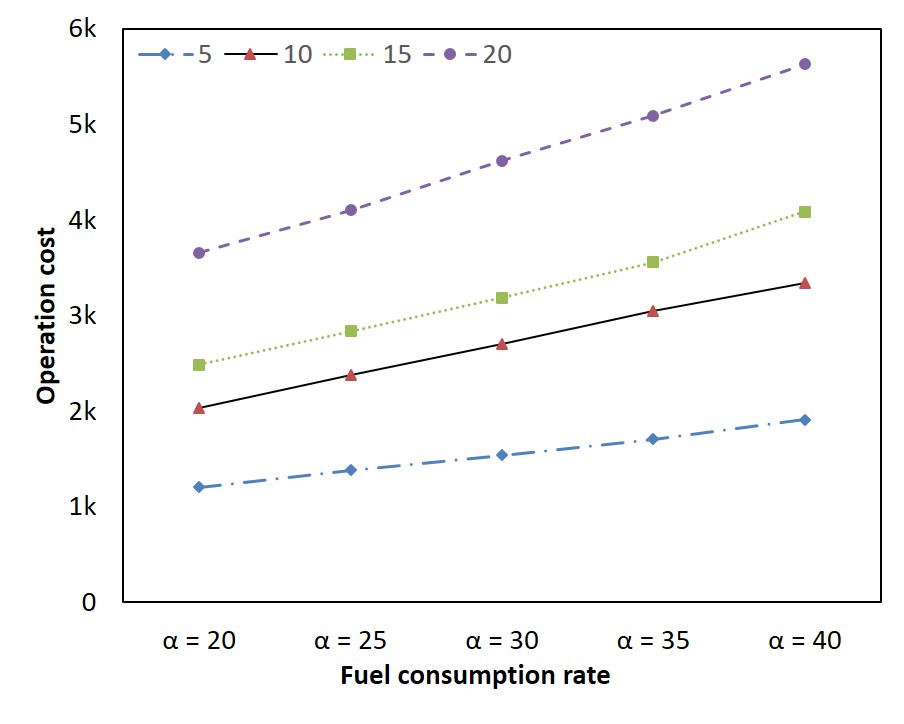}}\quad
	\subfloat[Platooning benefit in percentage]{\label{fig:fuelb}\includegraphics[width=0.48\textwidth]{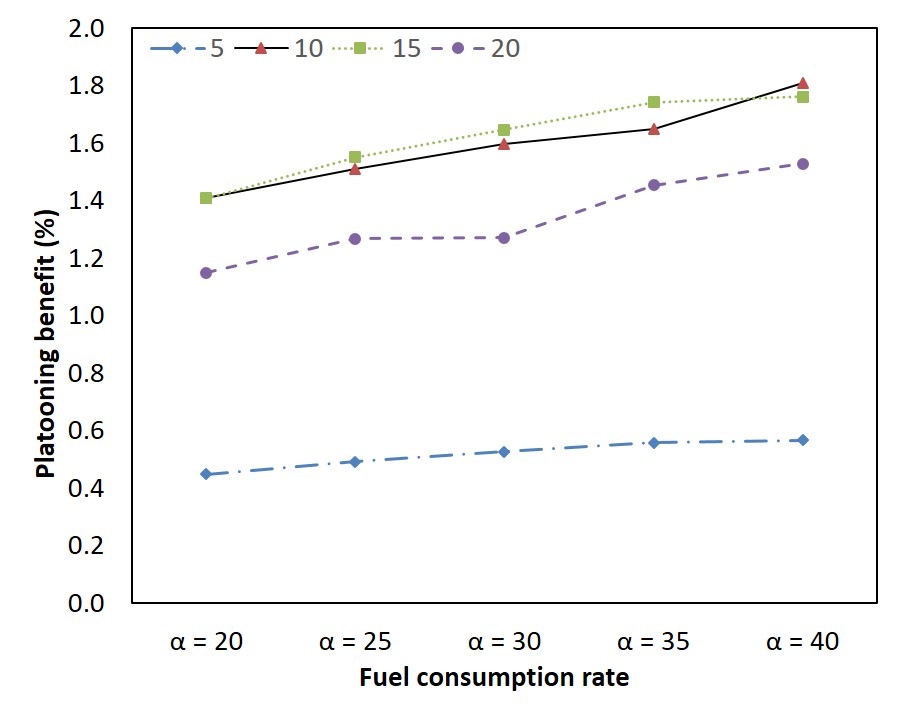}}\quad\\
	\caption{Influence of fuel consumption rate on operation cost and platooning benefit}
    \label{fig::fuel}
\end{figure}

In our previous experiments, we fixed the base fuel consumption rate, $\alpha$, at 30.7 by assuming that all trucks considered are typical Class 8 trucks. However, the fuel consumption rate may vary for different kinds of trucks, different weathers, and road conditions. We also anticipate witnessing a drop in energy consumption rate in dollars when electric trucks become more prevalent. Therefore, in this subsection, we capture the influence of fuel consumption rate on the system performance and platooning benefit. 
We chose 20, 25, 30, 35, and 40 as our tested levels for the fuel consumption rate, while other parameters remained the same across all test cases. Just as we illustrated before, for each test case with a specific fuel consumption rate and problem size, we took the average of 5 scenarios. The results in operation cost and platooning benefit are therefore obtained and visually presented in Figure \ref{fig::fuel}.

We can observe the relationship between the fuel consumption rate and the operation cost from Figure \ref{fig:fuela}. The curves of the first three cases (5, 10, and 15 customer nodes) show that the operation cost is positively linearly related to the fuel consumption rate. In the meanwhile, for the curves of platooning benefit shown in Figure \ref{fig:fuelb}, some turbulence happens when the fuel consumption rate is at 30 for the 20-customer case and at 35 for the 10-customer case, but an overall linear relationship can still be observed between the platooning benefit and the fuel consumption rates. In addition, the increasing trend tends to decline as the fuel consumption rate rises for most cases.  For instance, the increase in platooning benefit when the fuel consumption rate rises from 35 to 40 is smaller than the increase when the rate rises from 20 to 25. The shares of platooning benefit are capped at a certain threshold that is lower than 1.8\% in this example. 

Explaining that the percentages are calculated by the platooning benefit divided by the total operation cost. The 1.8\% platooning benefit in percentage may sound minimal, but the number is still significant considering the large value of the operation cost. In addition, we set the platoon saving ratio to be 10\%, and thus the percentage of platooning benefit can at most be 7.5\% when dispatch cost is omitted. The 7.5\% can only be obtained when each truck travels exactly the same route under the platooning and the non-platoon case, and trucks form a maximum-size platoon, which is selected as 4, on every link. Therefore, 1.8\% appears to be a reasonable number for the platooning benefit in real instances, and platooning benefit can still play an essential role in energy savings and carbon emission reduction.

\subsubsection{Unit vehicle dispatch cost}
In the previous experiments, we used $271$ dollars per truck to represent the fixed unit vehicle dispatch cost, and the value is selected based on the average salary of a truck driver in the United States. However, as the salary level varies across different countries and even different companies, the unit vehicle dispatch cost will also be different in different places. Therefore, in this subsection, we want to capture the influences of the unit vehicle dispatch cost on operation cost and platooning benefit. 

\begin{figure}[htbp]
	\centering
	\subfloat[Operation cost]{\label{fig:disa}\includegraphics[width=0.48\textwidth]{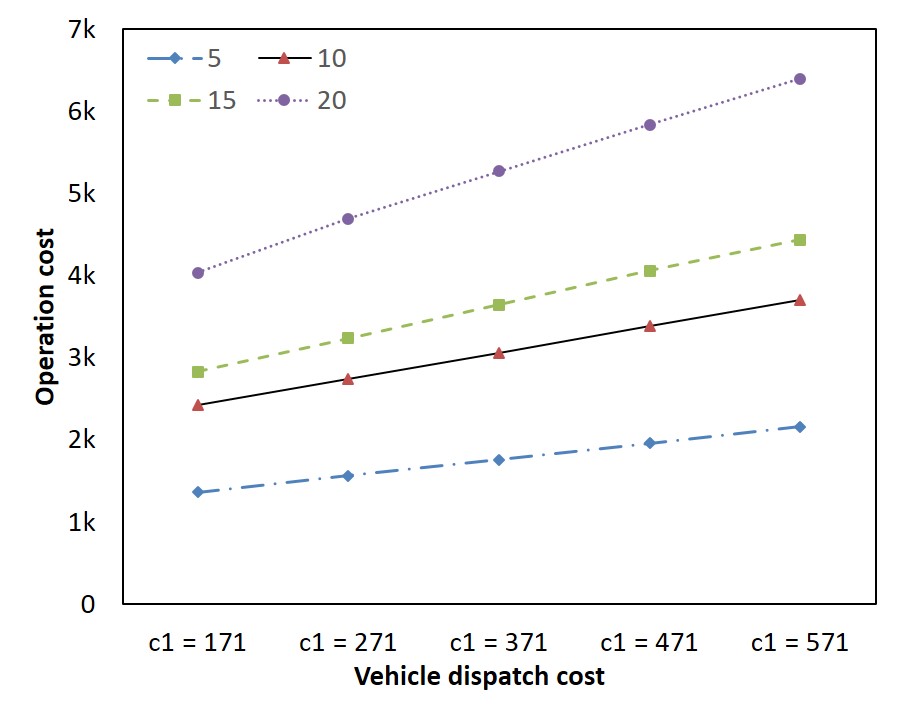}}\quad
	\subfloat[Platooning benefit in percentage]{\label{fig:disb}\includegraphics[width=0.48\textwidth]{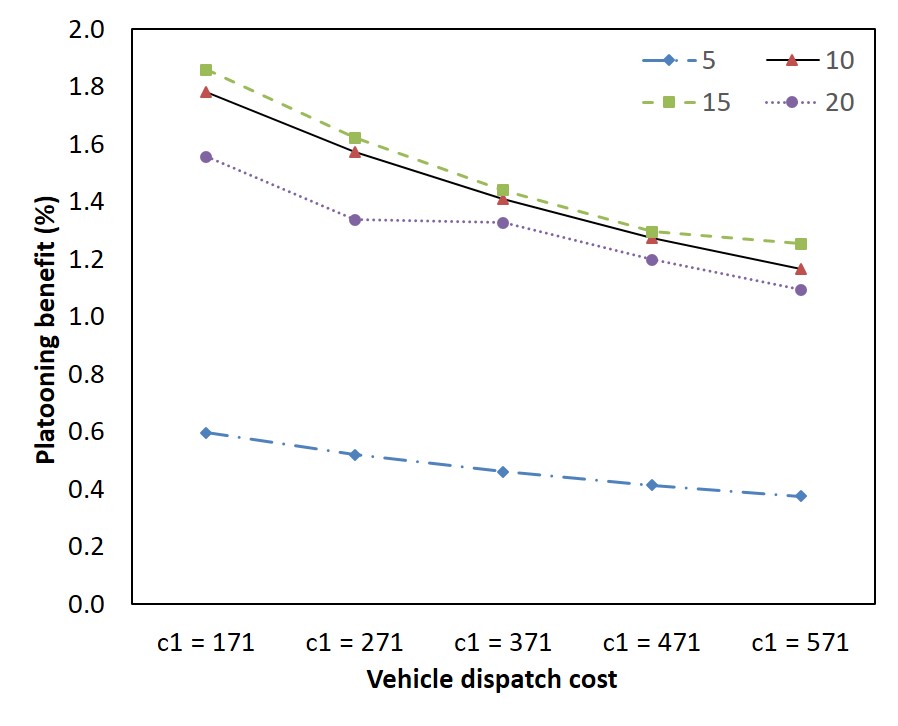}}\quad\\
	\caption{Influence of unit vehicle dispatch cost on operation cost and platooning benefit}
    \label{fig::dispatch}
\end{figure}

We chose the costs from five different levels, which are $171, 221, 271, 321$, and $371$, to cover the salary levels of most cases. We used the same test cases as the previous experiment based on the Yangtze River Delta network, but the difference was that we modified the parameter $c_1$ instead of $\alpha$ this time, and we fixed $\alpha$ as the default setting, which was 30.7. We tested the same case on different $c_1$ values to make their results comparable. The results of operation cost and platooning benefit are visually shown in Figure \ref{fig::dispatch}. We also tabulate the number of dispatched trucks in Table \ref{tab:dispveh}.

\begin{table}[htbp]
\caption{Average number of dispatched trucks under different cases}
\label{tab:dispveh}
\centering
\begin{threeparttable}
\begin{tabular}{cccccc}
\toprule

\textbf{Number of customer nodes}  & \textbf{$c_1 = 171$} & \textbf{$c_1 = 271$} & \textbf{$c_1 = 371$}  & \textbf{$c_1 = 471$} & \textbf{$c_1 = 571$}\\ 
\hline
\midrule
5 & 2 & 2 & 2 & 2 & 2\\
10 & 3.2& 3.2 & 3.2 & 3.2& 3.2\\
15 & 4.2 & 4.2 & 4.2 & 4.2 & 3.8\\
20 & 6.6 & 6.6 & 5.6 & 5.6 & 5.6\\
\hline

\bottomrule
\end{tabular}
\end{threeparttable}
\end{table}

According to Figure \ref{fig::dispatch}, we can observe the increase in unit dispatch cost greatly enlarges the total operation cost. For instance, as the unit vehicle dispatch cost increases from 171 to 571, the total operation cost of the 20-customer case rises from 4k to 6.5k, which is a more than 50\% increase. Combining with the results tabulated in Table \ref{tab:dispveh}, we can use simple math to further calculate that the total dispatch cost when $c_1 = 571$ occupies approximately $5.6 \times \dfrac{571}{6393} = 50.0\%$ of the total operation cost, telling that the dispatch cost occupies a huge portion of the total operation cost to some extent. In addition, as the average number of dispatched trucks remains the same for the 5-customer and 10-customer cases, we can claim that the increment in the total operation cost mainly comes from the increase in vehicle dispatch cost. Such a phenomenon possibly tells that the dispatching and routing strategy of trucks have not altered when the changes in unit dispatch cost are not significant and the customer demand size is limited. The reasons for such a phenomenon include that the central coordinator is nearly impossible to use fewer trucks to serve all the demands and that the savings from dispatching fewer trucks cannot compensate for the detours. However, the number of dispatched trucks does decrease when the dispatch cost increases from 471 to 571 for the 15-customer case, and from 271 to 371 for the 20-customer case. Visually speaking, a non-smooth disturbance can be witnessed at the transitions from $c_1 = 271$ to $c_1 = 371$, and $c_1 = 471$ to $c_1 = 571$, for the 20 and 15-customer cases, respectively, which reflects the alternation in dispatching strategy to maintain the platooning benefit. This implies that when the dispatch cost is too high, the central coordinator will tend to alter its dispatching strategy, such as using fewer drivers with more detours to save money. We can also witness that the relationship between the operation cost and the unit dispatch cost is positively linear, and the slope becomes steeper as the problem size scales. Though fewer trucks are dispatched, and more platooning benefit is generated, we can still witness an approximately linear increase in the total operation cost possibly because of the resulting detours by using fewer trucks.  Last but not least, the results shown in Figure \ref{fig:disb} do provide us a glance that the platoon savings in percentage can be more noteworthy for some cases, exceeding 1.8\% for the case of 15 and when the dispatch cost is 171,

\subsubsection{Platoon size limit}
\begin{figure}[htbp]
	\centering
	\subfloat[Operation cost]{\label{fig:psizea}\includegraphics[scale=0.5]{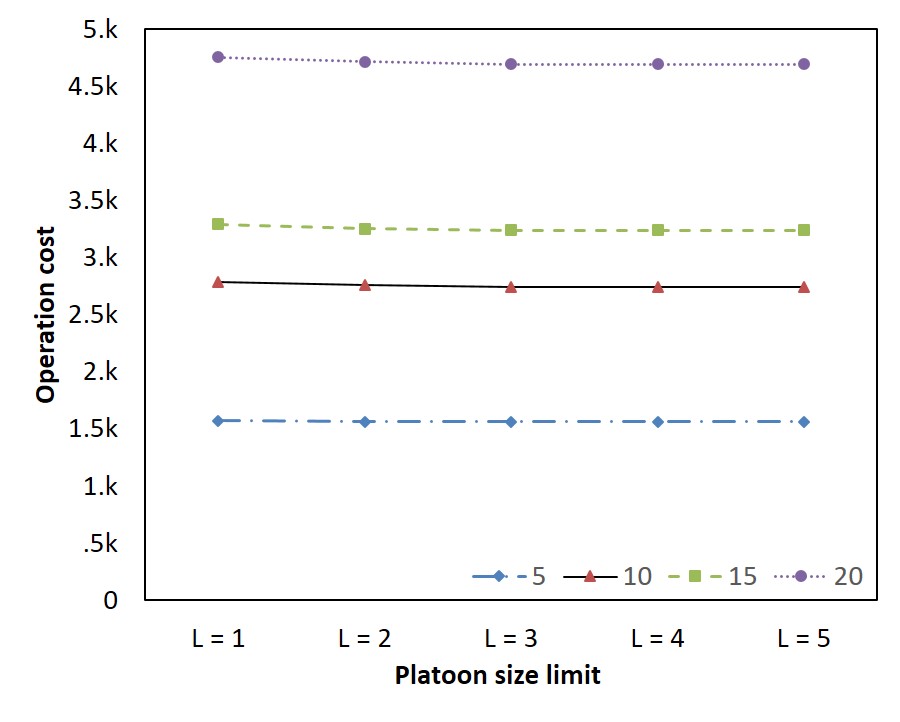}}\quad
	\subfloat[platooning benefit]{\label{fig:psizeb}\includegraphics[scale=0.5]{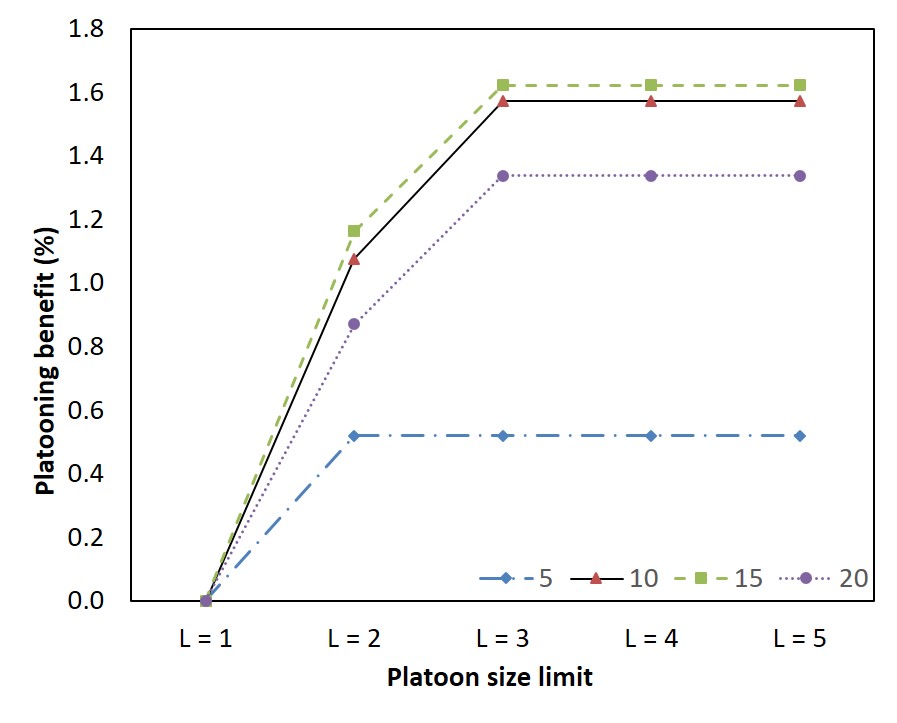}}\quad\\
	\caption{Influence of platoon size limit on operation cost and platooning benefit}
    \label{fig::platoonsize}
\end{figure}

In the previous experiments, we set the platoon size limit, $L$, to be 4 as recommended by \cite{laborcost}. In this subsection, we aim to investigate how will the planning strategy alter when different platoon sizes are applied. Therefore, we test 4 cases with 5 different platoon size limits, which are choosing from $[1,2,3,4,5]$. To be mentioned, when the platoon size limit equals 1, then it is equivalent to saying that no platoon formation is allowed. Therefore, such a case can represent the scenario without the platooning feature. All the results are visually presented in Figure \ref{fig::platoonsize}.

The results coincide with our expectations as the platooning benefit is none for all cases when the platoon size limit is set to 1, which indicates the scenario without platooning. Based on the curves shown in Figure \ref{fig:psizeb}, the platooning benefit is non-decreasing as the platoon size increases. As Figure \ref{fig:psizeb} illustrates, 2-truck platoons are formed for the 5-customer case as the platooning benefit occurs at the "$L = 2$" node and remains the same for the rest of the curve. 3-truck platoons are formed for the other cases as the platooning benefit step-wisely increases from $L=1$ to $L=2$, and from $L=2$ to $L = 3$, and the platooning benefit stagnates after the platoon size reaches 3. To summarize, the platoon size limit does affect the platooning formations and benefits, but it does not necessarily mean that the savings will definitely increase as the platoon size limit elevates because we can witness that the savings remain unchanged as the platoon size limit exceeds 3 in our experiments. Such a phenomenon infers that there exists a soft upper limit for the platoon savings and the longest platoon that can be formed, no matter how large the platoon size limit, especially when the demand size is not large. However, allowing longer platoons will generally bring savings in the operation cost for carriers, and the savings will possibly be more significant when the demand size scales because more trucks will be dispatched and more platoons may be formed.

\subsubsection{Truck capacity}

\begin{figure}[htbp]
	\centering
	\subfloat[Operation cost]{\label{fig:capa}\includegraphics[scale=0.5]{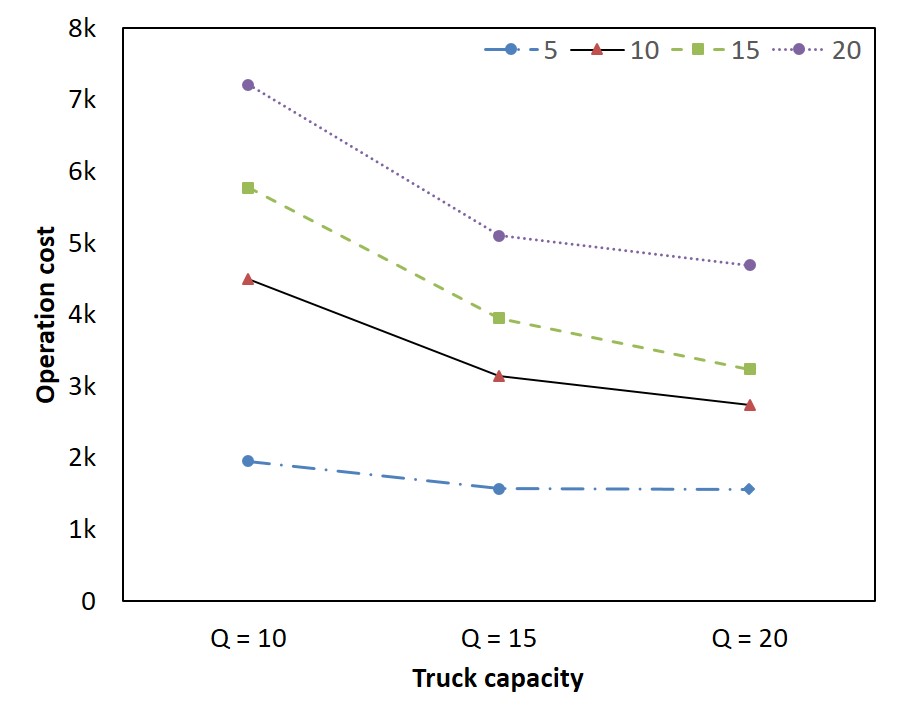}}\quad
	\subfloat[platooning benefit in percentage]{\label{fig:capb}\includegraphics[scale=0.5]{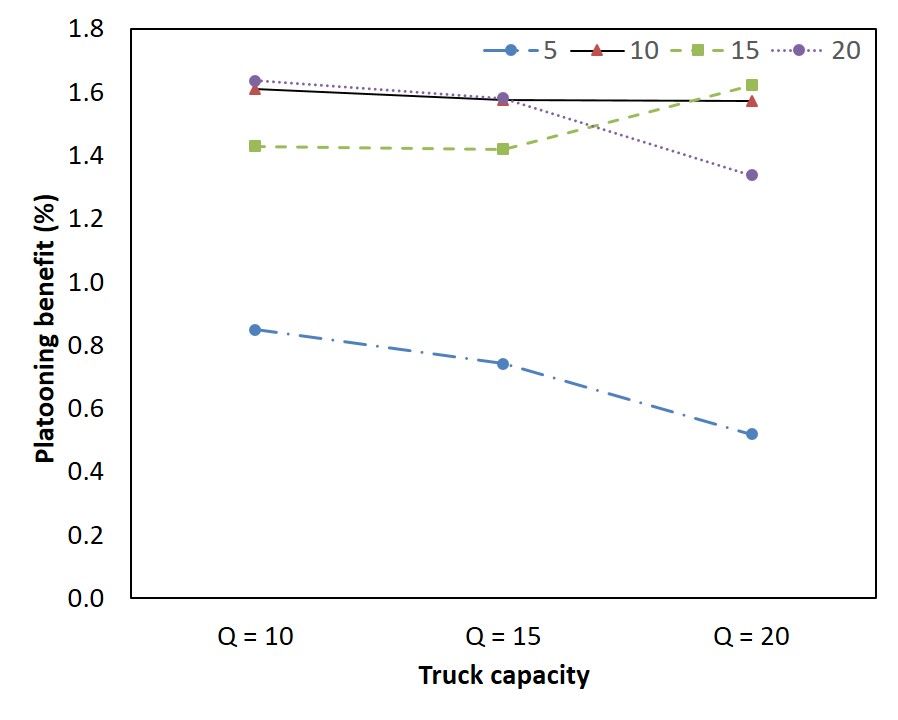}}\quad
        \subfloat[Absolute platooning benefit]{\label{fig:capc}\includegraphics[scale=0.5]{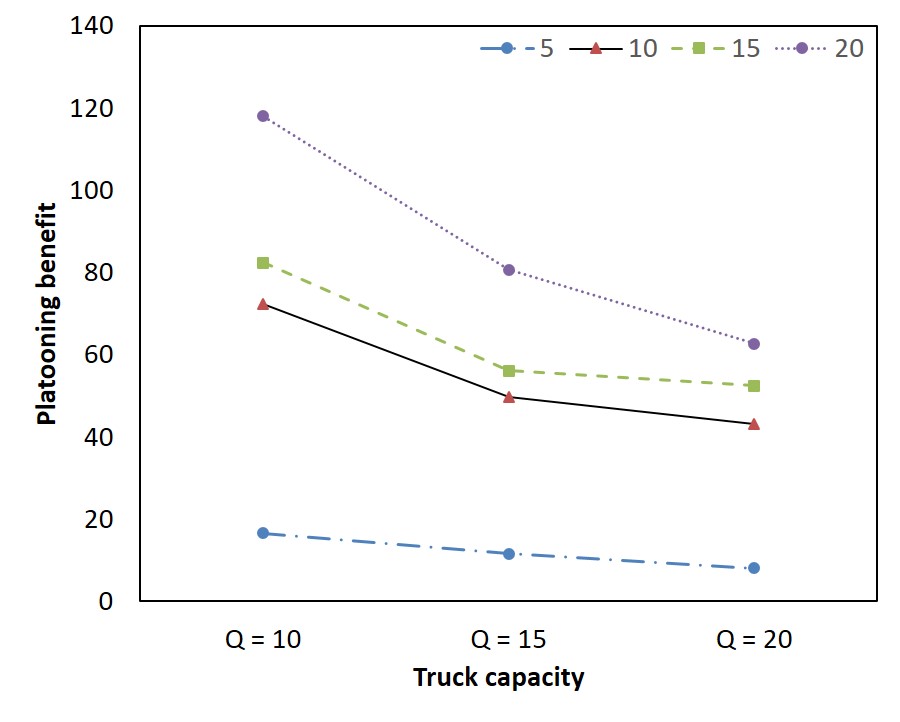}}\quad\\
	\caption{Influence of truck capacity on operation cost and platooning benefit}
    \label{fig::capacity}
\end{figure}

Assuming that we were dealing with typical 30-ton trucks (fully loaded) in real instances, similar to Class 8 trucks in the U.S., we fixed the truck's static weight at 10 and its capacity at 20. However, different types of trucks are used in different logistics service companies, and thus the truck capacity may vary among various scenarios. Therefore, in this subsection, we aim to check the impact of truck capacity on system performance and platooning benefit. Without loss of generality, we selected 10, 15, and 20 to be the three different capacity levels. Fixing the other parameters, the results in operation cost and platooning benefit under these three different capacity levels are visually presented in Figure \ref{fig::capacity}.

According to Figure \ref{fig:capa}, the operation cost exponentially rises as the truck capacity shrinks, as the differences between the results under $Q = 10$ and $Q = 15$ are more remarkable than the differences between the results under $Q = 15$ and $Q = 20$. Such exponential growth is universally true for all cases, and it is more significant for larger cases. One possible explanation for such a phenomenon is that the percent decrease in capacity is non-linear, 25\% from $Q = 20$ to $Q = 15$, and 33\% from $Q = 15$ to $Q = 10$. Another possible explanation is that the size of customer demands is randomly chosen within the range from 1 to 10, so setting the truck capacity to 10 may require significantly more trucks to fulfill the same amount of demands 
. Such an explanation can be further proved by the average number of dispatched trucks which is tabulated in Table \ref{tab:capveh}: the number of dispatched trucks increases drastically as the capacity diminishes, except for the case of 5 customer nodes, possibly because two trucks, even with capacity of 10, are enough to serve all the 5 customer nodes in our generated cases as the problem size is small.

\begin{table}[htbp]
\caption{Average number of dispatched trucks with different capacities}
\label{tab:capveh}
\centering
\begin{threeparttable}
\begin{tabular}{cccc}
\toprule

\textbf{Number of customer nodes}  & \textbf{$Q = 10$} & \textbf{$Q = 15$} & \textbf{$Q = 20$}  \\ 
\hline
\midrule
5 & 2 & 2 & 2 \\
10 & 7& 4.8 & 3.2 \\
15 & 9.2 & 7.2 & 4.2 \\
20 & 12.6 & 9.4 & 6.6 \\
\hline

\bottomrule
\end{tabular}
\end{threeparttable}
\end{table}

Evaluating the impact of truck capacity on platooning benefit, we may anticipate the platooning benefit to decrease as truck capacity enlarges because of fewer trucks dispatched, resulting in less platooning potential. According to Figure \ref{fig:capb}, the curves of 5, 10, and 20-customer cases coincide with our expectation, but the curve of the 15-customer case behaves differently. In this case, the platooning benefit in percentage increases as the truck capacity enlarges. To better understand the cause of this abnormal increase, we further visually present the absolute values of average platooning benefit in Figure \ref{fig:capc}, and it is shown that the absolute values of platooning benefit do strictly decrease as truck capacity enlarges. Therefore, we can infer that the abnormal increase in the platooning benefit is due to the huge difference in dispatched trucks. Since the number of dispatched trucks decreases from 7.2 to 4.2 when the capacity enlarges from 15 to 20, which is approximately a $\dfrac{7.2-4.2}{7.2} = 41.7\%$ decrease, such a huge change in dispatched trucks resulted in much more savings in dispatch cost and the total operation cost, thus underscoring the percentage of platooning benefit even if slight drop in absolute value has been witnessed. 

\subsubsection{Platoon saving ratio}
\begin{figure}[htbp]
	\centering
	\subfloat[Operation cost]{\label{fig:psavea}\includegraphics[scale=0.5]{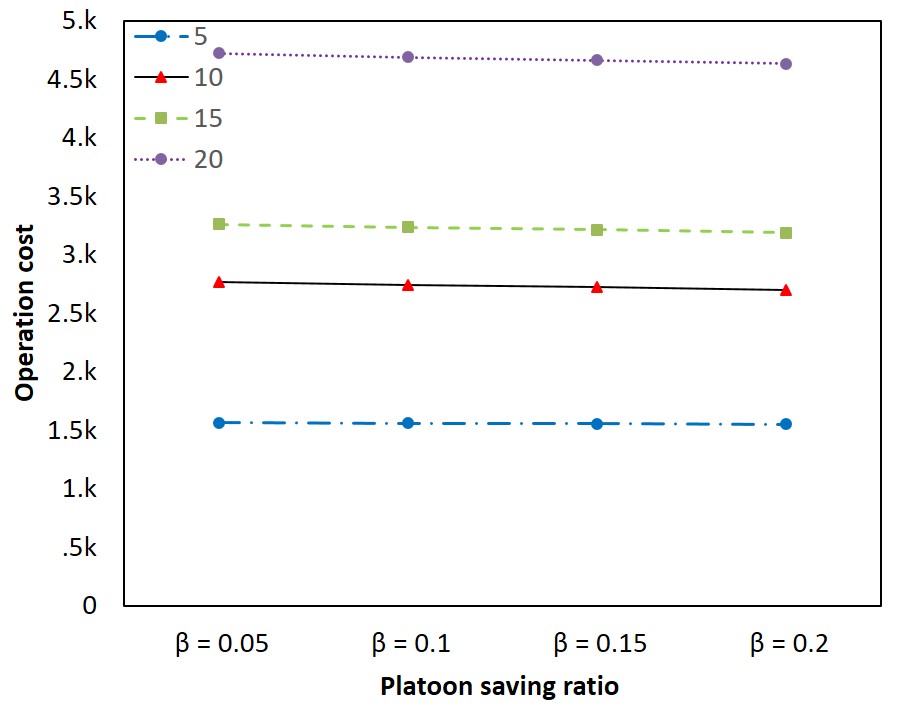}}\quad
	\subfloat[platooning benefit]{\label{fig:psaveb}\includegraphics[scale=0.5]{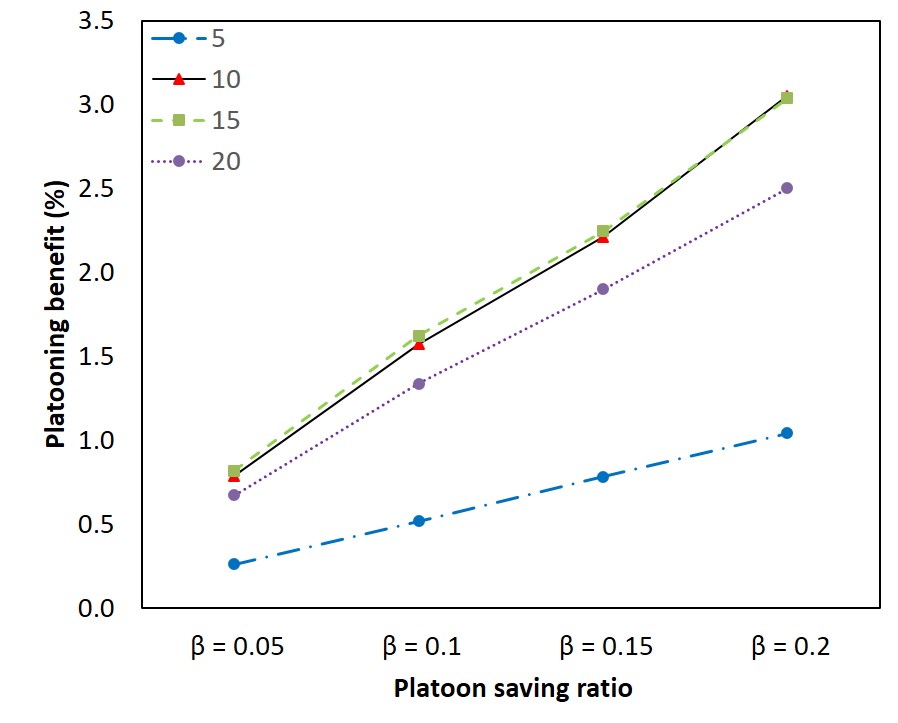}}\quad\\
	\caption{Influence of platoon saving ratio on operation cost and platooning benefit}
    \label{fig::platoonsave}
\end{figure}

The platoon saving ratio is often chosen from the range of 0.05 to 0.15 in most truck platooning studies, and 0.1 is the most frequently used one. We also realize that the platoon saving factor may vary under different road circumstances and truck categories, and it can be further elevated by closer headways between trucks due to the advancement of connected vehicles. Therefore, it is worth investigating the impact of the platoon saving ratio on operation cost and platooning benefit. We set up four platoon saving levels, which are 0.05, 0.1, 0.15, and 0.2 in the experiments. The results of these four cases are visually presented in Figure \ref{fig::platoonsave}.

Given the curves shown in Figure \ref{fig:psavea}, the influence of the platoon saving ratio on the total operation cost is relatively marginal because its increase does not yield a significant decline in the operation cost. However, while the platoon saving ratio increases, it drastically levels up the proportion of platooning benefit, as shown in Figure \ref{fig:psaveb}. For instance, the platooning benefit in percentage nearly doubles as the platoon saving ratio doubles. Furthermore, given the curves in Figure \ref{fig:psaveb}, we witness that the relation between the platooning benefit and the platooning saving ratio is approximately linear. Based on these findings, we can infer that the marginal decrease in total operation cost is mainly contributed by the additional platooning benefit, because, based on calculation, the marginal decrease in total operation cost is roughly the same as the increase in platooning benefit. 

\subsubsection{Time window tolerance}

\begin{figure}[htbp]
	\centering
	\subfloat[Operation cost]{\label{fig:twa}\includegraphics[scale=0.5]{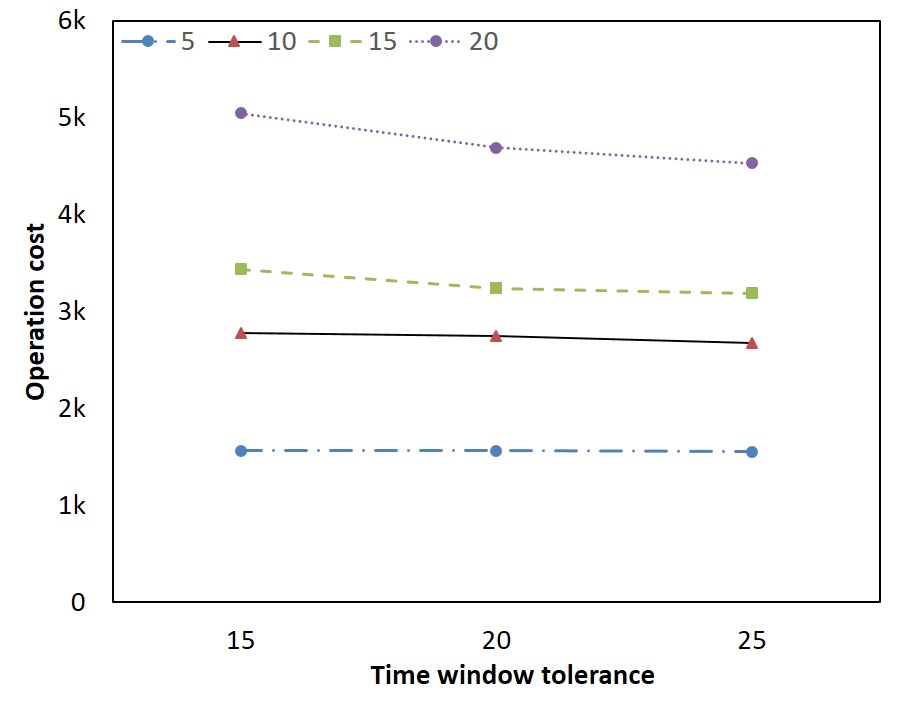}}\quad
	\subfloat[platooning benefit]{\label{fig:twb}\includegraphics[scale=0.5]{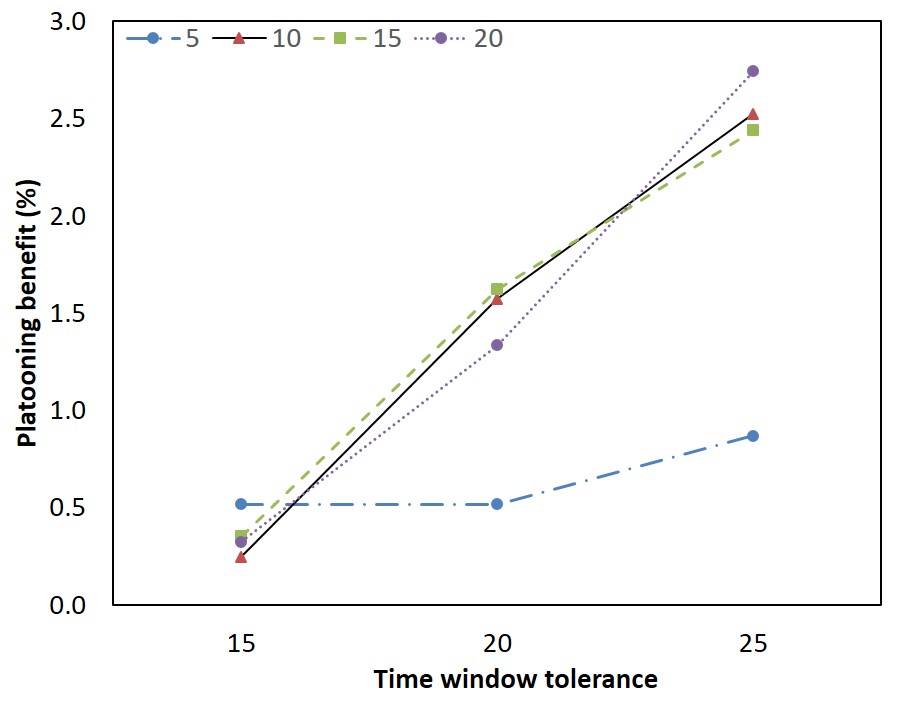}}\quad\\
	\caption{Influence of time window size on operation cost and platooning benefit}
    \label{fig::timewindow}
\end{figure}

In previous sections, we set the time window size between the latest departure time and the earliest arrival of each customer node to be at least 20, representing 20 hours. We can ensure the problem's feasibility with such a setup while maintaining a decent amount of platoon potential. In this subsection, we plan to alter the time window sizes to investigate their impact on platoon potential, which is reflected by the resulting value of the platooning benefit. We set a total of three levels, which are 15, 20, and 25. Therefore, during the time window generation of each customer, the latest departure time will be at least 15, 20, or 25 larger than its earliest arrival time. As a result, the additional time of 15, 20, and 25 can be understood as the time to travel with other customer demands together or the time to form platoons. Hence, we define this kind of additional time as time window tolerance in our experiments (because they are not the least necessary time to maintain problem feasibility). With different time window tolerances and other parameters fixed, the obtained results in operation cost and platooning benefit are provided in Figure \ref{fig::timewindow}.

Given the downtrend displayed in Figure \ref{fig:twb}, we can witness that the platooning benefit diminishes significantly as the time window tolerance shrinks. Comparing the curve of the 20-customer case and the 5-customer case, we can find that the percentage decrease in platooning benefit for the 20-customer case is about 2.5\% when the time window tolerance decreases from 25 to 15. However, the decrease in platooning benefit is even smaller than 0.5\% for the 5-customer case. Therefore, we can further highlight that the time window tolerance is much more crucial for cases with more customers and that platoon potential is also immense. On the other hand, a smaller time window tolerance yields an increase in the total operation cost, proved by the curves shown in Figure \ref{fig:twa}, but the increase is at a maximum of about 10\% even for the largest 20-customer case, which is limited. Therefore, the remarkable shrinking in platooning benefit in percentage is not mainly caused by the slight increase in the total operation cost but by the decrease in the absolute value of platooning benefit. Consequently, we can summarize that imposing a tighter time window tolerance will eliminate a considerable amount of platooning potential, thus leading to a decline in platooning benefit and an alternation of dispatching and routing strategies as more trucks may travel alone. 

\section{Conclusion}\label{sec:con}
This study discusses the planning of truck platooning for a capacitated road-network vehicle routing problem and highlights the benefits of applying the platooning feature for such a multi-delivery problem. A MIP model, defined as RCVRPTW-TP-R is proposed to optimize the dispatch and routing plan of trucks in a way to minimize the total operation cost of the central coordinator or a carrier company. The operation cost includes vehicles' dispatch cost and their on-road energy consumption costs, which are also influenced by trucks' weight and carrying loads. Specifically, the plan generated by the model will identify the serving truck for each customer node, the number of dispatched trucks, and the trucks' itineraries with time schedules so as to serve the customers within the required time windows. We propose a 3-stage solution algorithm embedded with a dynamic programming approach, a Modified Insertion Heuristic, and a "route-then-schedule" scheme to derive a high-quality solution with a reasonable amount of time. We further conduct numerical experiments on a virtual small network and a real-world network, i.e., the Yangtze River Delta network, to validate the feasibility of our model and demonstrate the platooning benefit. Comparisons between CPLEX and our proposed algorithm are presented to highlight the performance and effectiveness of our solution approach, and the results of most cases show that our proposed approach outperforms CPLEX both in time and solution quality. Sensitivity analyses have also been presented to quantify the impact of several key parameters on the operation plans and costs. It is revealed that the fuel consumption rate and the platoon size limit significantly influence the platooning benefit, but the impact of the platoon size limit on the total operation cost is limited. Unit dispatch cost influences the platooning benefit negatively, and it may alter the routing and dispatching strategy, thus inducing remarkable changes in the operation cost. In addition, platooning benefit is found to be positively correlated with time window tolerance and platoon saving ratio, but the impact of these factor on total operation cost are relatively small. Truck capacity has also proved to be important as its increase can effectively reduce the operation cost.

This study can be further extended or enhanced in several directions. From the perspective of problem formulation, goods transfer among trucks can be taken into consideration in formulating the problem to better depict the delivery service in the real-world scenario. Goods transfer is a valuable point to be further investigated because transshipments are universal in logistics service, and it can greatly improve logistics efficiency and save delivery costs, as mentioned by \cite{WOLFINGER2021105110}, \cite{sync1}, and \cite{syncbb}. In addition, goods transfer matches our basic model because we have already taken the truckload and weight into consideration, and thus involving goods transfer will not be a complicated add-on to our model. However, the inclusion of these factors may further complicate the problem and thus a more effective solution approach may be required. To develop a possibly more efficient solution algorithm, a more advanced feedback mechanism or perturbation strategy can be derived to further contract the solution space and fasten the solving process. For instance, we only shuffle the customer grouping results when the scheduling problem is infeasible, but we can further apply the shuffling mechanism on the feedback loop when the scheduling problem is feasible to broaden the solution space and avoid trapping at a local optimum.

\bibliographystyle{apalike} 
\footnotesize\bibliography{main}
\end{document}